\documentclass[12pt]{article}
\usepackage{mathrsfs}
\usepackage{epic,eepic,epsf,epsfig}
\usepackage{amsfonts,srcltx,mathrsfs}
\usepackage{amsmath}
\usepackage{amssymb}
\usepackage{amsbsy}
\usepackage{float}
\usepackage{graphicx}
\usepackage{amsfonts}
\usepackage{color}

\usepackage{url}
\newtheorem{lem}{Lemma}[section]
\newtheorem{theorem}[lem]{Theorem}

\newtheorem{cor}[lem]{Corollary}

\newtheorem{prob}[lem]{Problem}

\newtheorem{case}{Case}
\newtheorem{prop}[lem]{Proposition}
\newtheorem{rem}[lem]{Remark}

\def\a{\alpha}  \def\g{\gamma} \def\d{\delta} 
\def\r{\rho} \def\s{\sigma}

\def\olg{\overline G}

 \def\lg{\langle} \def\rg{\rangle}
\def\nd{\mathrel{\bigm|\kern-.7em/}}

\def\f{\noindent}

 \def\G{\hbox{\rm G}}

\def\Aut{\hbox{\rm Aut}}

\def\Aut{\hbox{\rm Aut}}

\def\Core{\hbox{\rm Core}}

\def\demo{\f {\bf Proof.}\hskip10pt}

\newcommand{\qed}{\mbox{\raisebox{0.7ex}{\fbox{}}} \vspace{4truemm}}
\def\mz{{\mathbb Z}}

\begin{document}

\title{Existence of regular $3$-polytopes of order $2^n$}

\author{ \\ Dong-Dong Hou$^{a}$, Yan-Quan Feng*$^{a}$, Dimitri Leemans$^{b}$\\
$^{a}${\small Department of Mathematics, Beijing Jiaotong University, Beijing,
100044, P.R. China}\\
$^{b}${\small {D\'epartement de Math\'ematique, Universit\'e Libre de Bruxelles, 1050 Bruxelles Belgium}}}

\date{}
\maketitle

\footnotetext{*Corresponding author.E-mails: yqfeng$@$bjtu.edu.cn, 16118416$@$bjtu.edu.cn,dleemans$@$ulb.ac.be
}
\begin{abstract}

In this paper, we prove  that for any positive integers $n, s, t$ such that $n \geq 10$, $s, t \geq 2$ and $n-1 \geq s+t$, there exists a regular polytope with Schl\"afli type $\{2^s, 2^t\}$ and its automorphism group is of order $2^n$. Furthermore, we classify regular polytopes with automorphism groups of order $2^n$ and Schl\"afli types $\{4, 2^{n-3}\}, \{4, 2^{n-4}\}$ and $\{4, 2^{n-5}\}$, therefore giving a partial answer to a problem proposed by Schulte and Weiss in [Problems on polytopes, their groups, and realizations, Periodica Math. Hungarica 53(2006) 231-255].

\bigskip
\f {\bf Keywords:} Regular polytope, $2$-group, automorphism group.
\end{abstract}

\section{Introduction  }

Classifications of abstract regular polytopes have been a subject of interest for several decades.
One path has been to fix (families of) groups of automorphisms and determine the abstract regular polytopes having these groups as full automorphism groups.
Some striking results have been obtained, for instance for the symmetric and alternating groups. Fernandes and Leemans classified abstract regular polytopes of rank $n-1$ and $n-2$ for $S_n$~\cite{fl,sympolcorr} and more recently, they extended this classification to rank $n-3$ and $n-4$ with Mixer~\cite{flm}.
Cameron, Fernandes, Leemans and Mixer showed that the highest rank of an abstract regular polytope with full automorphism group an alternating group $A_n$ is $\lfloor (n-1)/2 \rfloor$ when $n\geq 12$~\cite{CFLM2017}, and thanks to two previous papers of Fernandes, Leemans and Mixer~\cite{flm1,flm2}, this bound is known to be sharp.
More recently,  Gomi, Loyola and De Las Pe\~{n}as determined the non-degenerate string C-groups of order $1024$ in~\cite{sc1024}.

There exists a well known one-to-one correspondence between abstract regular polytopes and string C-groups. We therefore work with string C-groups as it is more convenient and easier to define them than abstract regular polytopes. In this paper, we study $2$-groups acting on regular polytopes. The starting point of our research was the following problem proposed by Schulte and Weiss in~\cite{Problem}.

\begin{prob}
Characterize the groups of orders $2^n$ or $2^np$, with $n$ a positive integer and $p$ an odd prime, which are automorphism groups of regular or chiral polytopes?
\end{prob}

Conder~\cite{SmallestPolytopes} showed that if $\mathcal{P}$ is a regular $3$-polytope with Schl\"afli type $\{k_1,k_2\}$, then $|\Aut(\mathcal{P})| \geq 2k_1k_2$. If $\mathcal{P}$ has Schl\"afli type $\{2^s,2^t\}$ and $|\Aut(\mathcal{P})|=2^n$, then $n-1\geq s+t$. In this paper, we first show the following theorem.

\begin{theorem}\label{existmaintheorem}
For any positive integers $n, s, t$ such that $n \geq 10$, $s, t \geq 2$ and $n-1 \geq s+t$, there exists a string C-group of order $2^n$ with Schl\"afli type $\{2^s, 2^t\}$.
\end{theorem}

Cunningham and Pellicer~\cite{GD2016} classified the  regular $3$-polytopes $\mathcal{P}$ for the case when $|\Aut(\mathcal{P})|=2k_1k_2$. Note that if $|\Aut(\mathcal{P})|=2^n$ and $k_1=4$ then $k_2 \leq 2^{n-3}$. As a special case, Cunningham and Pellicer~\cite{GD2016} obtained the classification of regular $3$-polytopes with automorphism groups of order $2^n$ and Schl\"afli type $\{4, 2^{n-3}\}$, and this was also given in Loyola~\cite{Loyola} by using the classification of $2$-groups with a cyclic subgroup of order $2^{n-3}$~\cite{Classfictionpgroup}. We prove the result again independently by using new techniques that are described in this paper, and the techniques work well for more classifications. In particular, we further
classify the regular $3$-polytopes with automorphism groups of order $2^n$ and Schl\"afli types $\{4, 2^{n-4}\}$ and $\{4, 2^{n-5}\}$ in this paper.

To state the next result, we need to define some groups:
\begin{small}
\begin{itemize}\setlength{\parskip}{0pt}
  \item [$G_1$]$=\lg \r_0, \r_1, \r_2 \ |\ \r_0^2, \r_1^2, \r_2^2, (\r_0\r_1)^{2^{2}}, (\r_1\r_2)^{2^{n-3}}, (\r_0\r_2)^2, [(\r_0\r_1)^2,\r_2] \rg$,
  \item [$G_2$]$=\lg \r_0, \r_1, \r_2 \ |\ \r_0^2, \r_1^2, \r_2^2, (\r_0\r_1)^{2^{2}}, (\r_1\r_2)^{2^{n-3}}, (\r_0\r_2)^2, [(\r_0\r_1)^2,\r_2](\r_1\r_2)^{2^{n-4}} \rg$,
  \item [$G_3$]$=\lg \r_0, \r_1, \r_2 \ |\ \r_0^2, \r_1^2, \r_2^2, (\r_0\r_1)^{2^{2}}, (\r_1\r_2)^{2^{n-4}}, (\r_0\r_2)^2, [(\r_0\r_1)^2, (\r_1\r_2)^2]\rg$,
  \item [$G_4$]$=\lg \r_0, \r_1, \r_2 \ |\ \r_0^2, \r_1^2, \r_2^2, (\r_0\r_1)^{2^{2}}, (\r_1\r_2)^{2^{n-4}}, (\r_0\r_2)^2, [(\r_0\r_1)^2, (\r_1\r_2)^2](\r_1\r_2)^{2^{n-5}} \rg$,
  \item [$G_5$]$=\lg \r_0, \r_1, \r_2    \ |\ \r_0^2, \r_1^2, \r_2^2, (\r_0\r_1)^{2^{2}}, (\r_1\r_2)^{2^{n-5}}, (\r_0\r_2)^2, [\r_0, (\r_1\r_2)^2]^2, [\r_0, (\r_1\r_2)^4]\rg$,
  \item [$G_6$]$=\lg \r_0, \r_1, \r_2 \ |\ \r_0^2, \r_1^2, \r_2^2, (\r_0\r_1)^{2^{2}}, (\r_1\r_2)^{2^{n-5}}, (\r_0\r_2)^2, [\r_0, (\r_1\r_2)^2]^2 (\r_1\r_2)^{2^{n-6}}, [\r_0, (\r_1\r_2)^4]\rg,$
  \item [$G_7$]$=\lg \r_0, \r_1, \r_2 \ |\ \r_0^2, \r_1^2, \r_2^2, (\r_0\r_1)^{2^{2}}, (\r_1\r_2)^{2^{n-5}}, (\r_0\r_2)^2, [\r_0, (\r_1\r_2)^2]^2, [\r_0, (\r_1\r_2)^4](\r_1\r_2)^{2^{n-6}}\rg,$
  \item [$G_8$]$=\lg \r_0, \r_1, \r_2 \ |\ \r_0^2, \r_1^2, \r_2^2, (\r_0\r_1)^{2^{2}}, (\r_1\r_2)^{2^{n-5}}, (\r_0\r_2)^2, [\r_0, (\r_1\r_2)^2]^2(\r_1\r_2)^{2^{n-6}},$ \\ \mbox{}\hskip 2.3cm $[\r_0, (\r_1\r_2)^4)](\r_1\r_2)^{2^{n-6}}\rg.$
\end{itemize}
\end{small}

\begin{theorem}\label{maintheorem}
For $n \geq 10$, let $\Gamma:=(G,\{\r_0,\r_1,\r_2\})$ be a string C-group of  order $2^n$. Then
\begin{enumerate}
\item[(1)] $\Gamma$ has type  $\{4,2^{n-3}\}$ if and only if $G \cong G_1$ or $G_2$;
\item[(2)] $\Gamma$ has type  $\{4,2^{n-4}\}$ if and only if $G \cong G_3$ or $G_4$;
\item[(3)] $\Gamma$ has type  $\{4,2^{n-5}\}$ if and only if $G \cong G_5, G_6, G_7$ or $G_8$.
\end{enumerate}
\end{theorem}

Let $n<10$. By \cite{atles1} or \cite{atles}, there is a unique string C-group of order $2^n$ with type $\{4,4\}$, and Theorem~\ref{maintheorem} is true for the types $\{4,2^{n-s}\}$ with $n-s\geq 3$ and  $s=3,4$ or $5$, except for the cases when $n=8$ or $9$ with $s=5$. For $n=8$ with $s=5$, there are four string C-groups with type $\{4,8\}$: two are $G_5$ and $G_6$, and the other two are $\lg \r_0, \r_1, \r_2 \ |\ \r_0^2, \r_1^2, \r_2^2, (\r_0\r_1)^{2^{2}}, (\r_1\r_2)^{2^{8}}, (\r_0\r_2)^2,  [(\r_0\r_1)^2, (\r_1\r_2)^2](\r_1\r_2)^4 \rg$ and $\lg \r_0, \r_1, \r_2 \ |\ \r_0^2$, $\r_1^2, \r_2^2, (\r_0\r_1)^{2^{2}}, (\r_1\r_2)^{2^{3}}, (\r_0\r_2)^2, [((\r_1\r_2)^2)^{\r_0}, \r_1\r_2](\r_1\r_2)^4\rg$. For $n=9$ with $s=5$, there are six string C-groups with type $\{4,16\}$: four are $G_i$ with $5\leq i\leq 8$, and the other two are
$\lg \r_0, \r_1, \r_2    \ |\ \r_0^2, \r_1^2, \r_2^2, (\r_0\r_1)^{2^{2}}, (\r_1\r_2)^{2^{4}}, (\r_0\r_2)^2,   [(\r_0\r_1)^2, (\r_1\r_2)^2](\r_1\r_2)^4 \rg$ and
$\lg \r_0, \r_1, \r_2 \ |\ \r_0^2$, $\r_1^2, \r_2^2, (\r_0\r_1)^{2^{2}}, (\r_1\r_2)^{2^{4}}, (\r_0\r_2)^2,    [(\r_0\r_1)^2, (\r_1\r_2)^2](\r_2\r_1)^4,
  [\r_1, \r_0, \r_2, \r_1, \r_0, \r_1, \r_0] \rg$.

\section{Background results}\label{backgroud}

\subsection{String C-groups}
Abstract regular polytopes and string C-groups are the same mathematical objects. The link between these objects may be found for instance in~\cite[Chapter 2]{ARP}.
We take here the viewpoint of string C-groups because it is the easiest and the most efficient one to define abstract regular polytopes.

Let $G$ be a group and let $S=\{\rho_0,\cdots,\rho_{d-1}\}$ be a generating set of involutions of $G$.
For $I\subseteq\{0,\cdots,d-1\}$, let $G_I$ denote the group generated by $\{\rho_i:i\in I\}$.
Suppose that
\begin{itemize}
\item[*]  for any $i,j\in \{0, \ldots, d-1\}$ with $|i-j|>1$, $\rho_i$ and $\rho_j$ commute (the \emph{string
property});
\item[*] for any $I,J\subseteq\{0,\cdots,d-1\}$,
$G_I\cap G_J=G_{I\cap J}\ \  (\mbox{the \emph{intersection property}})$.
\end{itemize}
Then the pair $(G,S)$ is called a {\em string C-group of rank} $d$ and the {\em order} of $(G,S)$ is simply the order of $G$. If $(G,S)$  only satisfies the string property, it is called a {\em string group generated by involutions} or \emph {sggi}. By the intersection property, $S$ is a minimal generating set of $G$.  It is known that string $C$-groups are the same thing as automorphism groups of regular polytopes~\cite[Section 2E]{ARP}. The following proposition is straightforward, and for details,  one may see \cite{MC}.

\begin{prop}\label{intersection}
The intersection property for a string $C$-group $(G,S)$ of rank $3$ is equivalent to that $S$ is a minimal generating set of $G$ and $\lg \r_0, \r_1\rg \cap \lg \r_1, \r_2\rg = \lg \r_1\rg $.
\end{prop}

The \emph{i-faces} of the regular d-polytope
associated with $(G,S)$ are the right cosets of the distinguished subgroup $G_i = \lg \r_j \ |\ j \neq i\rg$ for each $i=0,1,\cdots,d-1$, and two faces are incident just when they intersect as cosets. The {\em (Schl\"afli) type} of $(G,S)$ is the ordered set $\{p_1,\cdots, p_{d-1}\}$, where $p_i$ is the order of $\r_{i-1}\r_i$. In this paper we always assume that each $p_i$ is at least 3 for otherwise the generated group is a direct product of two smaller groups. If that happens, the string C-group (and the corresponding abstract regular polytope) is called {\em degenerate}. The following proposition is related to degenerate string C-groups of rank $3$.

\begin{prop}\label{degenerate}
For $t \geq 1$, let
\begin{itemize}\setlength{\parskip}{-3pt}
  \item [$L_1$]$=\lg \r_0, \r_1, \r_2 \ |\ \r_0^2, \r_1^2, \r_2^2, (\r_0\r_1)^{4}, (\r_1\r_2)^{2}, (\r_0\r_2)^2\rg$,
  \item [$L_2$]$=\lg \r_0, \r_1, \r_2 \ |\ \r_0^2, \r_1^2, \r_2^2, (\r_0\r_1)^{2}, (\r_1\r_2)^{2^t}, (\r_0\r_2)^2 \rg$,
  \item [$L_3$]$=\lg \r_0, \r_1, \r_2 \ |\ \r_0^2, \r_1^2, \r_2^2, (\r_0\r_1)^{2^t}, (\r_1\r_2)^{2}, (\r_0\r_2)^2 \rg$.
\end{itemize}
Then $|L_1|=16$, $|L_2|=|L_3|=2^{t+2}$. In particular,
the listed exponents are the true orders of the corresponding elements.
\end{prop}

The proof of Propostion~\ref{degenerate} is straightforward from the fact that  $L_2=\langle\r_0\rangle\times\langle\r_1,\r_2\rangle\cong \mz_2\times D_{2^{t+1}}$ and $L_3=\langle\r_0,\r_1\rangle\times\langle\r_2\rangle\cong D_{2^{t+1}}\times\mz_2$, where $D_{2^{t+1}}$ denotes the dihedral group of order $2^{t+1}$.

The following proposition is called the {\em quotient criterion} for a string C-group.

\begin{prop}{\rm \cite[Section 2E]{ARP}}\label{stringC}
Let $(\G,\{\r_0, \r_1, \r_2\})$ be an sggi, and let $\Lambda = (\lg \s_{0}, \s_{1}, \s_{2}\rg,$ $\{\s_0,\s_1,\s_2\})$  be a string C-group.
If the mapping $\r_j \mapsto \s_j$ for $j=0 ,1, 2$ induces a homomorphism $\pi : \G \rightarrow \Lambda$, which is one-to-one on the subgroup
 $\lg \r_0, \r_1\rg$ or on $\lg \r_1, \r_2\rg$, then $(\G, \{\r_0, \r_1, \r_2\})$ is also a string $C$-group.
\end{prop}

The following proposition gives some string C-groups with type $\{4, 4\}$, which is proved in \cite[Section 8.3]{HW} for $b\geq 2$ but it is also true for $b=1$ by {\sc Magma}~\cite{BCP97}.

\begin{prop}\label{type44}
For $b \geq 1$, let
\begin{itemize}\setlength{\parskip}{-3pt}
  \item [$M_1$]$=\lg \r_0, \r_1, \r_2 \ |\ \r_0^2, \r_1^2, \r_2^2, (\r_0\r_1)^{4}, (\r_1\r_2)^{4}, (\r_0\r_2)^2, (\r_2\r_1\r_0)^{2b}\rg$,
  \item [$M_2$]$=\lg \r_0, \r_1, \r_2 \ |\ \r_0^2, \r_1^2, \r_2^2, (\r_0\r_1)^{4}, (\r_1\r_2)^{4}, (\r_0\r_2)^2, (\r_1\r_2\r_1\r_0)^b\rg$.
\end{itemize}
Then $|M_1|=16b^2$ and $|M_2|=8b^2$. In particular,
the listed exponents are the true orders of the corresponding elements.
\end{prop}

\subsection{Permutation representation graphs and CPR graphs}
In~\cite{Pel2008}, Daniel Pellicer introduced CPR-graphs to give a permutation representation of string C-groups
(CPR stands for $C$-group Permutation Representation). These graphs are also sometimes called permutation representation graphs.

Let $G$ be a group and $S:= \{\rho_0,\ldots, \rho_{d-1}\}$ be a generating set of involutions of $G$.
Let $\phi$ be an embedding of $G$ into the symmetric group $S_n$ for some $n$.  The {\em permutation representation graph} ${\cal G}$ of $G$ determined by $\phi$ is the multigraph with $n$ vertices, and with edge labels in the set $\{0,\ldots,d-1\}$, such that any two vertices $v,w$ are joined by an edge of label $j$ if and only if $(v)((\rho_j)\phi)=w$.

If $(G,S)$ is a string C-group, then the permutation representation graph defined above is called a {\em CPR-graph} by Pellicer.

\subsection{Group theory}

Let $G$ be a group. For $x,y\in G$, we use $[x,y]$ as an abbreviation for the
{\em commutator} $x^{-1}y^{-1}xy$ of $x$ and $y$, and $[H, K]$ for the subgroup generated by all commutators $[x, y]$ with $x \in H$ and $y \in K$, when $H$ and $K$ are subgroups of $G$.
The following  proposition is a basic property of commutators and its proof is straightforward.

\begin{prop}\label{commutator}
Let $G$ be a group. Then, for any $x, y, z \in G$, $[xy, z]=[x, z]^y[y, z]$ and $[x, yz]=[x, z][x, y]^z$.
\end{prop}

The {\em commutator (or derived)} subgroup $G'$ of a group $G$ is the subgroup generated by all commutators $[x, y]$ for any $x, y \in G$. With Proposition~\ref{commutator}, it is easy to prove that if $G$ is generated by a subset $M$, then $G'$ is generated by all conjugates in $G$ of elements $[x_i, x_j]$ with $x_i, x_j \in M$; see ~\cite[Hilfsatz \uppercase\expandafter{\romannumeral3}.1.11]{GroupBooks} for example.

\begin{prop}\label{Derived}
Let $G$ be a group, $M \subseteq G$ and $G = \lg M \rg$. Then $G' = \lg [x_i, x_j]^g\ |\ x_i, x_j \in M, g \in G\rg$.
\end{prop}

The {\em Frattini subgroup}, denoted by $\Phi(G)$, of a finite group $G$ is defined to be the intersection of all maximal subgroups of $G$. Let $G$ be a finite $p$-group for a prime $p$,  and set $\mho_1(G) = \lg g^p\ |\ g \in G\rg$. The following theorem is the well-known Burnside Basis Theorem.

\begin{theorem}{\rm ~\cite[Theorem 1.12]{GroupBookss}}\label{burnside}
Let $G$ be a $p$-group and $|G: \Phi(G)| = p^d$.
\begin{itemize}
\item [(1)] $G/\Phi(G) \cong \mz_p^d$. Moreover, if $N \lhd G$ and $G/N$ is elementary abelian, then $\Phi(G) \leq N$.

\item [(2)] Every minimal generating set of $G$ contains exactly $d$ elements.

\item [(3)] $\Phi(G) = G' \mho_1(G)$. In particular, if $p=2$, then $\Phi(G) = \mho_1(G)$.
\end{itemize}
\end{theorem}

By Theorem~\ref{burnside}(2), we have the following important result.

\begin{rem}
A string $2$-group has $C$-group representations in only one rank.
\end{rem}

The unique cardinality of all  minimal generating set of a $2$-group $G$ is called the {\em rank} of $G$,  and denoted by $d(G)$. This is quite different from almost simple groups where in most cases if a group has string $C$-group representations of maximal rank $d$, then it has string $C$-group representations of ranks from $3$ to $d$. The only known exception is the alternating group $A_{11}$~\cite{flm1}.

For a subgroup $H$ of a group $G$, the {\em core} $\Core_G(H)$ of $H$ in $G$ is the largest normal subgroup of $G$ contained in $H$. The following result is called {\em Lucchini's theorem}.

\begin{prop}{\rm \label{core}\cite[Theorem 2.20]{GroupBook}}
Let $A$ be a cyclic proper subgroup of a finite group $G$, and let $K = \Core_{G}(A)$. Then $|A : K| < |G : A|$, and in
particular, if $|A| \geq |G : A|$, then $K > 1$.
\end{prop}

\section{Proof of Theorem~\ref{existmaintheorem}}\label{Theorem1.2}

Let $n \geq 10$, $s, t \geq 2$ and $n-s-t \geq 1$. Set $R(\r_0, \r_1, \r_2)= \{\r_0^2, \r_1^2, \r_2^2, (\r_0\r_1)^{2^{s}}, (\r_1\r_2)^{2^{t}}, (\r_0\r_2)^2$, $[(\r_0\r_1)^4, \r_2], [\r_0,(\r_1\r_2)^4]\}$ and define
$$H=\left\{
\begin{array}{ll}
\lg \r_0, \r_1, \r_2 \ |\ R(\r_0, \r_1, \r_2), [(\r_0\r_1)^2, \r_2]^{2^{\frac{n-s-t-1}{2}}}\rg, & n-s-t\mbox{ odd }\\
\lg \r_0, \r_1, \r_2 \ |\ R(\r_0, \r_1, \r_2), [(\r_0\r_1)^2, (\r_1\r_2)^2]^{2^{\frac{n-s-t-2}{2}}}\rg, & n-s-t\mbox{ even. }
\end{array}
\right.$$
To prove Theorem~\ref{existmaintheorem}, we only need to show that $H$ is a string C-group of order $2^{n}$ with Schl\"afli type $\{2^s, 2^{t}\}$. For convenience, write $o(h)$ for the order of $h$ in $H$.

Note that $\r_0$ commutes with $(\r_1\r_2)^4$ because $[\r_0,(\r_1\r_2)^4]=1$. Since $\langle\r_1,\r_2\rangle$ is a dihedral group, we have $(\r_1\r_2)^{\r_1}=(\r_1\r_2)^{\r_2}=(\r_1\r_2)^{-1}$. It follows that $\lg (\r_1\r_2)^4\rg \unlhd H$. Similarly, $\lg (\r_0\r_1)^4\rg \unlhd H$ as $[(\r_0\r_1)^4,\r_2]=1$.

Let $L_2=\lg \r_0, \r_1, \r_2 \ |\ \r_0^2, \r_1^2, \r_2^2, (\r_0\r_1)^{2}, (\r_1\r_2)^{2^t}, (\r_0\r_2)^2 \rg$. Clearly, $\r_0$ commutes with both $\r_1$ and $\r_2$ in $L_2$, and hence $[\r_0,(\r_1\r_2)^4]=1$. It is easy to see that the generators $\r_0, \r_1, \r_2$ in $L_2$ satisfy all relations in $H$. This implies that $L_2$ is a homomorphic image of $H$. By Proposition~\ref{degenerate},
$\r_1\r_2$ has order $2^t$ in $L_2$, and hence has order $2^t$ in $H$. It follows that $|H|=o((\r_1\r_2)^4) \cdot |H/\lg (\r_1\r_2)^4\rg|=2^{t-2} \cdot |H/\lg (\r_1\r_2)^4\rg|$.

Let $L_3=\lg \r_0, \r_1, \r_2 \ |\ \r_0^2, \r_1^2, \r_2^2, (\r_0\r_1)^{2^s}, (\r_1\r_2)^{2}, (\r_0\r_2)^2 \rg$. The element $\r_2$ commutes with both $\r_0$ and $\r_1$ in $L_3$, and hence $[(\r_0\r_1)^4, \r_2]=1$. Since $\r_0\r_2=\r_2\r_0$, Proposition~\ref{commutator} implies $[(\r_0\r_1)^2, \r_2]=[\r_1\r_0\r_1, \r_2]=[\r_0, \r_1\r_2\r_1]^{\r_1}=[\r_0, (\r_1\r_2)^2]^{\r_2\r_1}$.
Hence $[(\r_0\r_1)^2, \r_2]=1$ in $L_3$. Therefore the generators $\r_0, \r_1, \r_2$ in $L_3$ satisfy all relations in $H$. By Proposition~\ref{degenerate},
$\r_0\r_1$ has order $2^s$ in $L_3$, and hence has order $2^s$ in $H$.  It follows that $|H|=2^{s-2} \cdot |H/\lg (\r_0\r_1)^4\rg|$.

To finish the proof of Theorem~\ref{existmaintheorem}, we are left to prove that $|H|=2^n$.

\medskip
\f {\bf Case 1:} $s=2$.
We distinguish two cases, namely the case where $n-t$ is odd and the case where $n-t$ is even.

Assume that $n-t$ is odd. Then $H=\lg \r_0, \r_1, \r_2 \ |\ R(\r_0, \r_1, \r_2), [(\r_0\r_1)^2, \r_2]^{2^{\frac{n-t-3}{2}}}\rg$. Since $\r_0\r_2=\r_2\r_0$, we have $[(\r_0\r_1)^2, \r_2]=(\r_1\r_0\r_1\r_2)^2$. It follows that $[(\r_0\r_1)^2, \r_2]^{2^{\frac{n-t-3}{2}}}=(\r_1\r_0\r_1\r_2)^{2^{\frac{n-t-1}{2}}}$.
Note that $(\r_0\r_1)^4=1$ and $\lg (\r_1\r_2)^4 \rg \unlhd H$. Thus $H/\lg (\r_1\r_2)^4\rg\cong H_1$, where $H_1=\lg \r_0, \r_1, \r_2 \ |\ \r_0^2, \r_1^2, \r_2^2, (\r_0\r_1)^{2^2},  (\r_1\r_2)^{2^2}, (\r_0\r_2)^2,  (\r_1\r_0\r_1\r_2)^{2^{\frac{n-t-1}{2}}}\rg$. By Proposition~\ref{type44}, $|H_1|=8\cdot (2^{\frac{n-t-1}{2}})^2=2^{n-t+2}$, and hence $|H|=2^{t-2} \cdot |H/\lg (\r_1\r_2)^4\rg|=2^{t-2}|H_1|=2^n$.

Assume that $n-t$ is even. Then $H=\lg \r_0, \r_1, \r_2 \ |\ R(\r_0, \r_1, \r_2), [(\r_0\r_1)^2, (\r_1\r_2)^2]^{2^{\frac{n-t-4}{2}}}\rg$. A similar argument as above gives rise to $H/\lg (\r_1\r_2)^4\rg\cong H_2=$ $\lg \r_0, \r_1, \r_2 \ |\ \r_0^2, \r_1^2$, $ \r_2^2, (\r_0\r_1)^{2^2},$ $  (\r_1\r_2)^{2^2}, (\r_0\r_2)^2,  [(\r_0\r_1)^2, (\r_1\r_2)^2]^{2^{\frac{n-t-4}{2}}}\rg$. Noting $(\r_0\r_1)^2=(\r_0\r_1)^{-2}$ and $(\r_1\r_2)^2=(\r_1\r_2)^{-2}$ in $H_2$, we have
$[(\r_0\r_1)^2, (\r_1\r_2)^2]^{2^{\frac{n-t-4}{2}}}=$ $(((\r_0\r_1)^2(\r_1\r_2)^2)^2)^{2^{\frac{n-t-4}{2}}}=$ $(((\r_0\r_1\r_2)^2)^2)^{2^{\frac{n-t-4}{2}}}=$ $(\r_0\r_1\r_2)^{2\cdot 2^{\frac{n-t-2}{2}}}$ because $\r_0\r_2=\r_2\r_0$, and so $H_2=\lg \r_0, \r_1, \r_2 \ |\ \r_0^2, \r_1^2$, $ \r_2^2, (\r_0\r_1)^{2^2},  (\r_1\r_2)^{2^2},$ $(\r_0\r_2)^2,  (\r_0\r_1\r_2)^{2\cdot 2^{\frac{n-t-4}{2}}}\rg$. By Proposition~\ref{type44}, $|H_2|=16 \cdot (2^{\frac{n-t-2}{2}})^2=2^{n-t+2}$, and $|H|=2^{t-2} \cdot |H/\lg (\r_1\r_2)^4\rg|=2^{t-2} \cdot |H_2|=2^n$.

\medskip
\f {\bf Case 2:} $s>2$.

Assume that $n-t-s$ is odd. Then $H=\lg \r_0, \r_1, \r_2 \ |\ R(\r_0, \r_1, \r_2), [(\r_0\r_1)^2, \r_2]^{2^{\frac{n-t-s-1}{2}}}\rg$. It follows $H/\lg (\r_0\r_1)^4\rg\cong H_3$, where $H_3=\lg \r_0, \r_1, \r_2 \ |\ \r_0^2, \r_1^2, \r_2^2, (\r_0\r_1)^{2^2},  (\r_1\r_2)^{2^t}, (\r_0\r_2)^2, $
$[\r_0, (\r_1\r_2)^4], [(\r_0\r_1)^2, \r_2]^{2^{\frac{(n-s+2)-t-3}{2}}}\rg$. By Case~1, $|H_3|=2^{n-s+2}$, and therefore $|H|=2^{s-2} \cdot |H/\lg (\r_0\r_1)^4\rg|=2^{s-2} \cdot |H_3|=2^n$.

Assume that $n-t-s$ is even. Then $H/\lg (\r_0\r_1)^4\rg\cong H_4$, where $H_4=\lg \r_0, \r_1, \r_2 \ |\ \r_0^2, \r_1^2, \r_2^2$, $(\r_0\r_1)^{2^2},  (\r_1\r_2)^{2^t}, (\r_0\r_2)^2, $
$[\r_0, (\r_1\r_2)^4], [(\r_0\r_1)^2, (\r_1\r_2)^2]^{2^{\frac{(n-s+2)-t-4}{2}}}\rg$. Then $|H_4|=2^{n-s+2}$ from Case~1, and therefore $|H|=2^{s-2} \cdot |H_4|=2^n$. \hfill\qed

\begin{cor} \label{cor3.1} The pairs $(G_1,\{\r_0,\r_1,\r_2\})$, $(G_3,\{\r_0,\r_1,\r_2\})$ and $(G_5,\{\r_0,\r_1,\r_2\})$, defined in Theorem~\ref{maintheorem}, are string $C$-groups of order $2^n$ with Schl\"afli type $\{4, 2^{n-3}\}$, $\{4, 2^{n-4}\}$ and $\{4, 2^{n-5}\}$, respectively.
\end{cor}

\demo By taking $(s, t)=(2, n-3),(2,n-4),(2,n-5)$ in the proof of Theorem~\ref{existmaintheorem}, we know that $(H_i,\{\r_0,\r_1,\r_2\})$ for $i=1,3,5$ are string $C$-groups of order $2^n$ with Schl\"afli type $\{4, 2^{n-3}\}$, $\{4, 2^{n-4}\}$ and $\{4, 2^{n-5}\}$ respectively, where

\begin{small}
\begin{itemize}\setlength{\parskip}{0pt}
  \item [$H_1$]$=\lg \r_0, \r_1, \r_2 \ |\ \r_0^2, \r_1^2, \r_2^2, (\r_0\r_1)^{2^{2}}, (\r_1\r_2)^{2^{n-3}}, (\r_0\r_2)^2, [\r_0,(\r_1\r_2)^4], [(\r_0\r_1)^2,\r_2] \rg$,
   \item [$H_3$]$=\lg \r_0, \r_1, \r_2 \ |\ \r_0^2, \r_1^2, \r_2^2, (\r_0\r_1)^{2^{2}}, (\r_1\r_2)^{2^{n-4}}, (\r_0\r_2)^2, [\r_0,(\r_1\r_2)^4], [(\r_0\r_1)^2, (\r_1\r_2)^2]\rg$,
    \item [$H_5$]$=\lg \r_0, \r_1, \r_2    \ |\ \r_0^2, \r_1^2, \r_2^2, (\r_0\r_1)^{2^{2}}, (\r_1\r_2)^{2^{n-5}}, (\r_0\r_2)^2, [\r_0, (\r_1\r_2)^4], [\r_0, (\r_1\r_2)^2]^2\rg$.
  \end{itemize}
\end{small}

Since $\r_0\r_2=\r_2\r_0$ and $((\r_0\r_1)^2)^{\r_0}=((\r_0\r_1)^2)^{\r_1}=(\r_0\r_1)^2$, by Proposition~\ref{commutator} we have the following identities in all $H_i$ and $G_i$ for $i=1,3,5$:
\begin{small}
\begin{itemize}\setlength{\parskip}{0pt}
  \item []$[(\r_0\r_1)^2, \r_2]=[\r_1\r_0\r_1, \r_2]=[\r_0, \r_1\r_2\r_1]^{\r_1}=[\r_0, (\r_1\r_2)^2]^{\r_2\r_1}$,
   \item []$[\r_0, (\r_1\r_2)^4]=[\r_0, (\r_1\r_2)^2][\r_0, (\r_1\r_2)^2]^{(\r_1\r_2)^2}=[(\r_0\r_1)^2, \r_2]^{\r_1\r_2}[(\r_0\r_1)^2, \r_2]^{(\r_1\r_2)^3}$,
    \item []$[(\r_0\r_1)^2,(\r_1\r_2)^2]=[(\r_0\r_1)^2,\r_2][(\r_0\r_1)^2,\r_1\r_2\r_1]^{\r_2}=
[(\r_0\r_1)^2,\r_2][(\r_0\r_1)^2,\r_2]^{\r_1\r_2}.$
  \end{itemize}
\end{small}

Clearly, $H_5=G_5$. In $G_1$,
$[\r_0, (\r_1\r_2)^4]=[(\r_0\r_1)^2, \r_2]^{\r_1\r_2}[(\r_0\r_1)^2, \r_2]^{(\r_1\r_2)^3}=1$ because $[(\r_0\r_1)^2,\r_2]=1$. Thus, $H_1=G_1$. To prove $H_3=G_3$, we only need to show that $[\r_0, (\r_1\r_2)^4]=1$ in $G_3$. Noting that $[(\r_0\r_1)^2, (\r_1\r_2)^2]=1$ in $G_3$, we have \begin{center}
$\begin{array}{l}
[(\r_0\r_1)^2, \r_2]^{\r_1}[(\r_0\r_1)^2, \r_2]=[\r_0,
(\r_1\r_2)^2]^{\r_2}[(\r_0\r_1)^2, \r_2]
=[\r_0, (\r_2\r_1)^2][(\r_0\r_1)^2, \r_2]\\
=\r_0(\r_1\r_2)^2\r_0(\r_2\r_1)^2\{(\r_0\r_1)^2\}\r_2(\r_0\r_1)^2\r_2
=\r_0\{(\r_0\r_1)^2\}(\r_1\r_2)^2\r_0(\r_2\r_1)^2\r_2(\r_0\r_1)^2\r_2
\\=\r_1\r_0\r_1(\r_1\r_2)^2\r_0\{\r_2(\r_1\r_2)^2\}(\r_0\r_1)^2\r_2
=\r_1\r_0\r_1(\r_1\r_2)^2\r_0\r_2(\r_0\r_1)^2(\r_1\r_2)^2\r_2
\\=\r_1\r_0\r_1\r_1\r_2\r_1\r_2\r_0\r_2\r_0\r_1\r_0\r_1\r_1\r_2\r_1\r_2\r_2=1,
\end{array}$
\end{center}
that is, $[(\r_0\r_1)^2, \r_2]^{\r_1}=[(\r_0\r_1)^2, \r_2]^{-1}$. On the other hand, since $1=[(\r_0\r_1)^2,\r_2\r_2]=[(\r_0\r_1)^2,\r_2][(\r_0\r_1)^2,\r_2]^{\r_2}$ and $1=[(\r_0\r_1)^2,(\r_1\r_2)^2]=
[(\r_0\r_1)^2,\r_2][(\r_0\r_1)^2,\r_2]^{\r_1\r_2}$, we have $[(\r_0\r_1)^2, \r_2]^{\r_1}=([(\r_0\r_1)^2, \r_2]^{\r_2})^{-1}=[(\r_0\r_1)^2, \r_2]$ . It follows that $[(\r_0\r_1)^2, \r_2]^{2}=1$ and
$[(\r_0\r_1)^2,\r_2]^{\r_1\r_2}=[(\r_0\r_1)^2,\r_2]^{-1}$, and so $[\r_0, (\r_1\r_2)^4]=[(\r_0\r_1)^2, \r_2]^{\r_1\r_2}[(\r_0\r_1)^2, \r_2]^{(\r_1\r_2)^3}=[(\r_0\r_1)^2,\r_2]^{-2}=1$, as required. \hfill\qed

\section{Proof of Theorem~\ref{maintheorem}}\label{Theorem1.3}

To prove Theorem~\ref{maintheorem}, we need the following lemmas.

\begin{lem}\label{quotient}
Let $(G,\{\r_0,\r_1,\r_2\})$ be a string C-group of type $\{2^s, 2^t\}$ with $2 \leq s \leq t$. Let $|G|=2^n$ and $2t \geq n-1$. Then $N = \lg (\r_1\r_2)^{2^{t-1}} \rg \unlhd G$ and $(\olg, \{\overline{\r_{0}}, \overline{\r_{1}}, \overline{\r_{2}}\})$ is a string $C$-group of type $\{2^s, 2^{t-1}\}$ and order $2^{n-1}$, where $\olg=G/N$ and $\overline{x}=xN$ for any $x\in G$.
\end{lem}

\demo Let $H = \lg \r_0\r_1, \r_1\r_2\rg$ be the rotation subgroup of $G$.
Then $|G:H| \leq 2$.
Since $\{\r_0,\r_1,\r_2\}$ is a minimal generating set of $G$, Theorem~\ref{burnside}~(2) implies $d(G)=3$, and since $H$ is generated by two elements, we have $|H| = 2^{n-1}$.

Let $M = \lg \r_1\r_2\rg$. Then $|M| = 2^t$ and $|M|^2\geq |H|$ as $2t \geq n-1$. Then Proposition~\ref{core} implies that $\Core _{H}(M)>1$. Since $M$ is cyclic and $|N| = 2$, $N$ is characteristic in $\Core _{H}(M)$, and so $\Core _{H}(M)\unlhd H$ implies $N \unlhd H$. Noting that $N$ lies in the center of the dihedral group $\lg \r_1, \r_2\rg$, we have $N^{\r_1} = N$, and hence $N \unlhd G$ because $G = \lg H, \r_1\rg$. Clearly, $\overline{G}=2^{n-1}$.

Since $t\geq 2$, we have $N\leq \mho_1(G) = \lg g^2 | g \in G\rg$, and by Theorem~\ref{burnside}, $N\leq \Phi(G)$ and $G/\Phi(G)\cong\mz_2^3$.
Thus, $G/\Phi(G)$ has rank $3$, and since $G/\Phi(G)\cong (G/N)/(\Phi(G)/N)$, $G/N$ has rank $3$, implying that $\{\overline{\r_0}, \overline{\r_1}, \overline{\r_2}\}$ is a minimal generating set of $\olg$. It follows that $\overline{\r_0}$, $\overline{\r_1}$, $\overline{\r_2}$ and $\overline{\r_0\r_2}$ are involutions. To prove that $\overline{G}$ has the intersection property, by Proposition~\ref{intersection} we only need to show $\lg \overline{\r_0}, \overline{\r_{1}} \rg \cap
\lg \overline{\r_{1}}, \overline{\r_{2}}\rg = \lg \overline{\r_{1}}\rg$.

Suppose $\lg \overline{\r_0}, \overline{\r_{1}} \rg \cap
 \lg \overline{\r_{1}}, \overline{\r_{2}}\rg \not= \lg \overline{\r_{1}}\rg$. Then there exist $x_1\in \lg \r_0,\r_1\rg$ and $x_2\in \lg \r_1,\r_2\rg$ such that $\overline{x_1} = \overline{x_2}\not\in \lg \overline{\r_1}\rg$, which implies $x_1\not\in
 \lg\r_1\rg$. Since $\lg \r_0, \r_1 \rg \cap
 \lg \r_1, \r_2\rg=\lg \r_1\rg$, we have $x_1 \neq x_2$ and  $x_1=x_2(\r_1\r_2)^{2^{t-1}}$ as $\overline{x_1} = \overline{x_2}$, which is impossible because otherwise $x_1=x_2(\r_1\r_2)^{2^{t-1}} \in \lg \r_0, \r_1\rg \cap \lg \r_1, \r_2\rg = \lg \r_1\rg$. Thus,  $\lg \overline{\r_0}, \overline{\r_{1}} \rg \cap
\lg \overline{\r_{1}}, \overline{\r_{2}}\rg = \lg \overline{\r_{1}}\rg$, as required.

To finish the proof, we are left to show $\overline{\r_0\r_1}$ and $\overline{\r_1\r_2}$ have order $2^s$ and $2^{t-1}$ respectively. Since $(G,\{\r_0,\r_1,\r_2\})$ has type $\{2^s,2^t\}$, $\r_0\r_1$ and $\r_1\r_2$ have order $2^s$ and $2^t$ respectively. Since $N\leq \langle \r_1\r_2\rangle$ and $|N|=2$, $\overline{\r_1\r_2}$ has order $2^{t-1}$ and $\overline{\r_0\r_1}$ has order $2^s$ or $2^{s-1}$.

Suppose  $\overline{\r_0\r_1}$ has order $2^{s-1}$. Then  $(\r_0\r_1)^{2^{s-1}}=(\r_1\r_2)^{2^{t-1}}\in \lg \r_0, \r_1\rg \cap \lg \r_1, \r_2\rg = \lg \r_1 \rg$, and hence $(\r_1\r_2)^{2^{t-1}}=\r_1$ because $\r_1\r_2$ has order $2^t$. It follows that $(\r_1\r_2)^{2^{t-1}-1}=\r_2$ and $(\r_1\r_2)^{2^t-2}=1$, a contradiction. Thus,  $\overline{\r_0\r_1}$ has order $2^s$. This completes the proof.
\hfill\qed

\begin{lem}\label{deviredgroup}
Let $G=\lg \r_0, \r_1, \r_2 \ |\ \r_0^2, \r_1^2, \r_2^2, (\r_0\r_2)^2\rg$. Then $G'=\lg [\r_0, \r_1], [\r_1, \r_2], [\r_0, \r_1]^{\r_2}\rg$.
\end{lem}

\demo Since $\r_0\r_2=\r_2\r_0$, Proposition~\ref{Derived} implies $G' = \lg [\r_0, \r_1]^{g}, [\r_2, \r_1]^{h} |\ g, h \in G\rg$. Since $\lg \r_0,\r_1\rg$ and $\lg \r_1,\r_2\rg$ are dihedral groups, we have $[\r_0, \r_1]^{\r_0}=((\r_0\r_1)^2)^{\r_0}=(\r_0\r_1)^{-2}=[\r_0, \r_1]^{-1}$, $[\r_0, \r_1]^{\r_1}=[\r_0, \r_1]^{-1}$, $[\r_1, \r_2]^{\r_1}=[\r_1, \r_2]^{-1}$ and $[\r_1, \r_2]^{\r_2}=[\r_1, \r_2]^{-1}$.

Set $L=\lg [\r_0, \r_1], [\r_1, \r_2], [\r_0, \r_1]^{\r_2}\rg$. Since $([\r_0, \r_1]^{\r_2})^{\r_0}=([\r_0, \r_1]^{\r_2})^{-1}$, $([\r_0, \r_1]^{\r_2})^{\r_2}=[\r_0, \r_1]$ and $([\r_0, \r_1]^{\r_2})^{\r_1}=\r_1\r_2\r_0\r_1\r_0\r_1\r_2\r_1=[\r_1, \r_2][\r_1, \r_0]^{\r_2}[\r_2, \r_1]$, we have $[\r_0, \r_1]^g \in L$ for any $g\in G$. Since $[\r_1, \r_2]^{\r_0} = \r_0\r_1\r_2\r_1\r_2\r_0=\r_0\r_1\r_0\r_1\r_1\r_2\r_1\r_2\r_2\r_1\r_0\r_1\r_0\r_2$
$=[\r_0, \r_1][\r_1, \r_2][\r_1, \r_0]^{\r_2}$, we have $[\r_1, \r_2]^h \in L$ for any $h\in G$. It follows that $G' \leq L$, and hence $G'=L$. \hfill\qed

\vskip 0.3cm
\f{\bf Proof of Theorem~\ref{maintheorem}(1):} For the sufficiency, we need to show that both $G_1$ and $G_2$ are string $C$-groups of order $2^n$, where $n\geq 10$ and

\vskip 0.1cm
\begin{small}
$\begin{array}{rl}
  G_1=& \lg \r_0, \r_1, \r_2 \ |\ \r_0^2, \r_1^2, \r_2^2, (\r_0\r_1)^{2^{2}}, (\r_1\r_2)^{2^{n-3}}, (\r_0\r_2)^2, [(\r_0\r_1)^2,\r_2] \rg,  \\
  G_2=& \lg \r_0, \r_1, \r_2 \ |\ \r_0^2, \r_1^2, \r_2^2, (\r_0\r_1)^{2^{2}}, (\r_1\r_2)^{2^{n-3}}, (\r_0\r_2)^2, [(\r_0\r_1)^2,\r_2](\r_1\r_2)^{2^{n-4}} \rg.
\end{array}$
\end{small}
\vskip 0.1cm

By Corollary~\ref{cor3.1}, $G_1$ is a string $C$-group of order $2^n$ and we are only left with $G_2$. However, to explain the method clearly, we prove the above fact again for $G_1$ using a permutation representation graph that is simple and easy to understand. Let $G=G_1$ or $G_2$. For convenience, we write $o(g)$ for the order of $g$ in $G$. We first prove the following claim.

\medskip
\f {\bf Claim:} $|G|\leq 2^n$.

Note that $G/G'$ is abelian and is generated by three involutions. Thus $|G/G'|\leq 2^3$. To prove the claim, it suffices to show $|G'|\leq 2^{n-3}$.

For $G=G_1$, we have $[(\r_0\r_1)^2, \r_2]=1$ and $\r_0\r_2=\r_2\r_0$, which implies $[\r_0, \r_1]^{\r_2}=[\r_0, \r_1][\r_1, \r_0][\r_0, \r_1]^{\r_2}=[\r_0, \r_1][(\r_0\r_1)^2, \r_2] = [\r_0,\r_1]$. Since $[\r_0,\r_1]^{\r_0}=[\r_0,\r_1]^{\r_1}=[\r_0,\r_1]^
{-1}$, we have $\lg [\r_0,\r_1]\rg\unlhd G$, and by Lemma~\ref{deviredgroup}, we have $G'=\lg [\r_0, \r_1], [\r_1, \r_2], [\r_0, \r_1]^{\r_2}\rg=\langle [\r_0, \r_1], [\r_1, \r_2]\rg=\lg [\r_0,\r_1]\rg\lg [\r_1,\r_2]\rg$. This implies that $|G'|\leq  |\lg [\r_0,\r_1]\rg||\lg [\r_1,\r_2]\rg|=o((\r_0\r_1)^2)o((\r_1\r_2)^2)\leq 2\cdot 2^{n-4}=2^{n-3}$, as required.

For $G=G_2$, we have $[(\r_0\r_1)^2,\r_2](\r_1\r_2)^{2^{n-4}}=1$ and $\r_0\r_2=\r_2\r_0$, which implies $[\r_0, \r_1]^{\r_2}=[\r_0, \r_1][(\r_0\r_1)^2, \r_2] = [\r_0, \r_1](\r_1\r_2)^{-2^{n-4}} \in \lg [\r_0, \r_1], [\r_1, \r_2]\rg$ and $[\r_1, \r_2]^{\r_0}=\r_1[(\r_0\r_1)^2, \r_2]\r_2\r_1\r_2=\r_1(\r_1\r_2)^{-2^{n-4}}\r_1(\r_1\r_2)^2=(\r_1\r_2)^{2^{n-4}}(\r_1\r_2)^2 \in \lg [\r_1, \r_2]\rg$. It follows that $\lg [\r_1, \r_2]\rg\unlhd G$ because $[\r_1,\r_2]^{\r_1}=[\r_1,\r_2]^{\r_2}=[\r_1,\r_2]^{-1}$, and by Lemma~\ref{deviredgroup}, $G'=\lg [\r_0, \r_1], [\r_1, \r_2], [\r_0, \r_1]^{\r_2}\rg=\langle [\r_0, \r_1],[\r_1, \r_2]\rg=\langle [\r_0, \r_1]\rg\lg [\r_1, \r_2]\rg$. In particular, $|G'|\leq  |\lg [\r_0,\r_1]\rg||\lg [\r_1,\r_2]\rg|=o((\r_0\r_1)^2)o((\r_1\r_2)^2)\leq 2\cdot 2^{n-4}=2^{n-3}$, as required.

Now we are ready to finish the sufficiency proof by considering two cases.
We use another method than the quotient method, based on permutation representation graphs. We give the details for $G_1$ as they are simpler than those of $G_2$ and might help the reader understand the case $G_2$.

\medskip
\f {\bf Case 1:} $G=G_1$.

\begin{figure}[H]
  \centering
  \includegraphics[width=9cm]{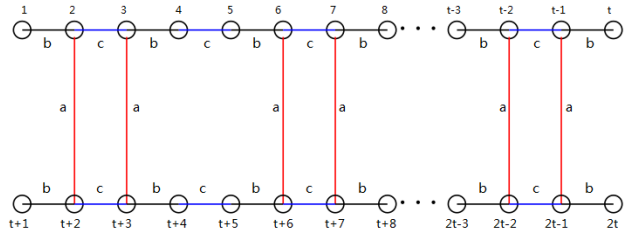}\\
  \caption{A permutation representation graph corresponding to $G_1$}\label{1}
\end{figure}

The key point is to construct a permutation group $A$ of order at least $2^n$ on a set $\Omega$ that is an epimorphic image of $G$, that is, $A$ has three generators, say $a, b, c$, satisfying the same relations as do $\r_0,\r_1,\r_2$. The permutation representation graph has vertex set $\Omega$ with $a$-, $b$- and $c$-edges. Recall that an $x$-edge ($x$=$a,b$ or $c$) connects two points in $\Omega$ if and only if $x$ interchanges them. It is easy to have such graphs when $n$ is small by taking $\Omega$ as the set of right cosets of the subgroup $\lg\r_0,\r_2\rg$ in $G$, where $\r_0$, $\r_1$ and $\r_2$ produce the $a$-, $b$- and $c$-edges, respectively. We give in Figure~\ref{1} a permutation representation graph for $G_1$ and explain below how it is constructed.

Set $t=2^{n-3}$ and write $i_{jt}^{k}=jt+4i+k$ where $0 \leq i \leq \frac{t}{4}-1, 0 \leq j \leq 1$ and $1\leq k \leq 4$. Then $a, b, c$ are permutations on the set $\{1,2,\cdots, 2^{n-2}\}$:

\vskip 0.2cm

\begin{small}
$\begin{array}{l}
a=\prod_{i=0}^{\frac{t}{4}-1}(i_{0}^{2},i_{t}^{2})(i_{0}^{3},i_{t}^{3}), ~~~\quad b=\prod_{i=0}^{\frac{t}{4}-1}(i_{0}^{1},i_{0}^{2})(i_{t}^{1},i_{t}^{2})(i_{0}^{3},i_{0}^{4})(i_{t}^{3},i_{t}^{4}),\\
c=(0_{0}^{1})(0_{t}^{1})((\frac{t}{4}-1)_{0}^{4})((\frac{t}{4}-1)_{t}^{4})\cdot\prod_{i=0}^{\frac{t}{4}-1}(i_{0}^{2},i_{0}^{3})(i_{t}^{2},i_{t}^{3})\cdot\prod_{i=0}^{\frac{t}{4}-2}(i_{0}^{4},(i+1)_{0}^{1})(i_{t}^{4},(i+1)_{t}^{1}).\\
\end{array}$
\end{small}
\vskip 0.2cm
Here, $(i+1)_{jt}^k=jt+4(i+1)+k$ for $0\leq i\leq \frac{t}{4}-2$. Note that $1$-cycles are also given in the product of distinct cycles of $c$ and this would be helpful to compute conjugations of some elements by $c$. It is easy to see that $a$ is fixed under conjugacy of $c$, that is, $a^c=a$. It follows $(ac)^2=1$. We further have

\vskip 0.2cm
\begin{small}
$\begin{array}{lcl}
  ab&=&\prod_{i=0}^{\frac{t}{4}-1}
  (i_{0}^{1},i_{0}^{2},i_{t}^{1},i_{t}^{2})(i_{0}^{3},i_{t}^{4},i_{t}^{3},i_{0}^{4}),\\

  bc&=&\prod_{i=0}^{1}(1+ti,3+ti,\cdots, t-1+ti, t+ti, t-2+ti, \cdots, 2+ti), \\

 (ab)^2&=&\prod_{i=0}^{\frac{t}{8}-1}
 (i_{0}^{1},i_{t}^{1})(i_{0}^{2},i_{t}^{2})(i_{0}^{3},i_{t}^{3})(i_{0}^{4},i_{t}^{4}).\\
\end{array}$
\end{small}
\vskip 0.2cm

Let $A=\langle a,b,c\rangle$. Clearly, $a^2=b^2=c^2=1$, $(ab)^4=1$ and $(bc)^{2^{n-3}}=1$. Furthermore, $(ab)^2$ is fixed under conjugacy of $c$, that is, $((ab)^2)^c=(ab)^2$, and hence $[(ab)^2,c]=1$. Clearly, $A$ is transitive on $\{1,2,\cdots, 2^{n-2}\}$, and the stabilizer $A_1$ has order at least $4$ because $a,c\in A_1$. This implies that $A$ is a permutation group of order at least $2^n$ and its generators $a,b,c$ satisfy the same relations as do $\r_0,\r_1,\r_2$ in $G$. Then there is an epimorphism $\phi:$ $G\mapsto A$ such that $\r_0^\phi=a$, $\r_1^\phi=b$ and $\r_2^\phi=c$. Since $|A|\geq 2^n$ and $|G|\leq 2^n$, $\phi$ is an isomorphism, implying $|G|=2^n$.

The generators $\r_0, \r_1, \r_2$ in $L_1:=\lg \r_0, \r_1, \r_2 \ |\ \r_0^2, \r_1^2, \r_2^2, (\r_0\r_1)^{4}, (\r_1\r_2)^{2}, (\r_0\r_2)^2\rg$ satisfy all relations in $G$. This implies that the map: $\r_0\mapsto\r_0$, $\r_1\mapsto\r_1$, $\r_2\mapsto\r_2$, induces a homomorphism from $G$ to $L_1$. By Proposition~\ref{degenerate}, $o(\r_0\r_1)=4$ in $L_1$ and hence $o(\r_0\r_1)=4$ in $G$, and by Proposition~\ref{stringC}, $(G,\{\r_0,\r_1,\r_2\})$ is a string C-group.

\medskip
\f {\bf Case 2:} $G=G_2$.

\begin{figure}[H]
  \centering
  \includegraphics[width=14cm]{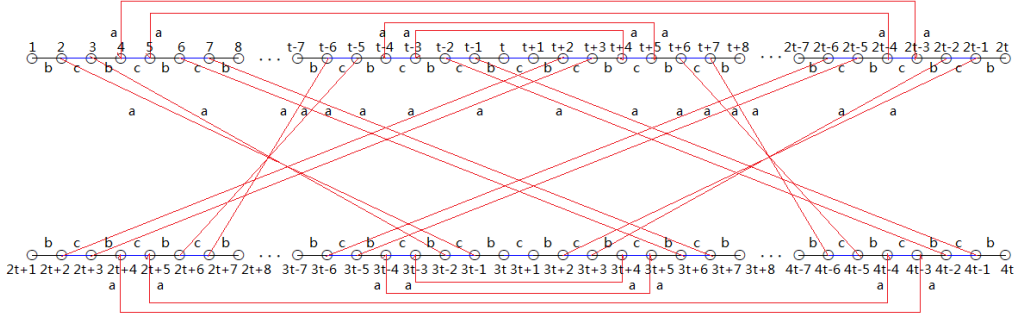}\\
  \caption{A permutation representation graph corresponding to $G_2$}\label{2}
\end{figure}

As in Case 1, we give in Figure~\ref{2} a permutation representation graph for $G_2$. Note that, as in Case~1, $bc$ consists of two paths of length $2t$ with alternating labels of $b$ and $c$, where $t=2^{n-4}$, and here all of the real complexity lies in the definition of $a$.

Write
$ci=\frac{t}{8}-i-1$ for $0\leq i\leq \frac{t}{8}-1$,
$i_{jt}^k=jt+8i+k$ and ${ci}_{jt}^k=jt+8ci+k$ for $0 \leq j \leq 3$ and $1 \leq k \leq 8$. Note that $0\leq i\leq \frac{t}{8}-1$ if and only if $0\leq ci\leq \frac{t}{8}-1$. Then $a, b, c$ are permutations on the set $\{1,2,\cdots, 2^{n-2}\}$:
\vskip 0.2cm
\begin{small}
$\begin{array}{rl}
  a=&\prod_{i=0}^{\frac{t}{8}-1}
    (i_t^2,i_{2t}^2)(i_0^2,ci_{2t}^7)(i_{3t}^2,ci_t^7)(i_0^7,i_{3t}^7)(i_t^3,i_{2t}^3)(i_0^3,ci_{2t}^6)
    (i_{3t}^3,ci_t^6)(i_0^6,i_{3t}^6)\\
    &(i_0^4,ci_t^5)(i_0^5,ci_t^4)(i_{2t}^4,ci_{3t}^5)(i_{2t}^5,ci_{3t}^4), \\

    b=&\prod_{j=0}^3\prod_{i=0}^{\frac{t}{8}-1}(i_{jt}^1,i_{jt}^2)(i_{jt}^3,i_{jt}^4)
  (i_{jt}^5,i_{jt}^6)(i_{jt}^7,i_{jt}^8),          \\

c=& \prod_{i=0}^{1} (0_{0+2ti}^{1}) ((\frac{t}{8}-1)_{t+2ti}^{8})((\frac{t}{8}-1)_{0+2ti}^{8},0_{t+2ti}^{1})
\prod_{j=0}^3(\prod_{i=0}^{\frac{t}{8}-1}(i_{jt}^2,i_{jt}^3)(i_{jt}^4,i_{jt}^5)
  (i_{jt}^6,i_{jt}^7) \cdot\\

& \prod_{i=0}^{\frac{t}{8}-2}(i_{jt}^8,(i+1)_{jt}^1)) .\\

\end{array}$
\end{small}\\
Here, $(i+1)_{jt}^1=jt+8(i+1)+1$ for $0\leq i\leq \frac{t}{8}-2$. It is easy to see that $a$ is fixed under conjugacy of $c$, that is, $a^c=a$. It follows $(ac)^2=1$.

Let $\a=a,b$, or $c$. Then $\a$ is an involution. Recall that $ci=\frac{t}{8}-i-1$. Since  $0\leq i\leq \frac{t}{8}-1$ if and only if $0\leq ci\leq \frac{t}{8}-1$, it is easy to see that if $\a$ interchanges $i_{j_1t}^{k_1}$ and $i_{j_2t}^{k_2}$ then $\a$ also interchanges $ci_{j_1t}^{k_1}$ and $ci_{j_2t}^{k_2}$, and if $\a$ interchanges $i_{j_1t}^{k_1}$ and $ci_{j_2t}^{k_2}$ then $\a$ also interchanges $ci_{j_1t}^{k_1}$ and $i_{j_2t}^{k_2}$. These facts are very helpful for the following computations.

\vskip 0.2cm
\begin{small}
$\begin{array}{lcl}
  ab&=&\prod_{i=0}^{\frac{t}{8}-1}
      (i_t^1,i_t^2,i_{2t}^1,i_{2t}^2)
      (i_0^1,i_0^2,{ci}_{2t}^8,{ci}_{2t}^7)
      (i_{3t}^1,i_{3t}^2,{ci}_t^8,{ci}_t^7)
      (i_0^7,i_{3t}^8,i_{3t}^7,i_0^8) \\

    & &(i_0^3,{ci}_{2t}^5,i_{3t}^3,{ci}_t^5)
      (i_0^5,{ci}_t^3,{ci}_{2t}^4,i_{3t}^6)
      (i_0^6,i_{3t}^5,{ci}_{2t}^3,{ci}_t^4)
      (i_0^4,{ci}_t^6,i_{3t}^4,{ci}_{2t}^6),\\

  bc&=&\prod_{i=0}^{1}(1+2ti,3+2ti,\cdots, 2t-1+2ti, 2t+2ti, 2t-2+2ti, \cdots, 2+2ti), \\

 (ab)^2&=&\prod_{i=0}^{\frac{t}{8}-1}

     (i_t^1,i_{2t}^1)(i_0^1,ci_{2t}^8)(i_t^8,ci_{3t}^1)(i_0^8,i_{3t}^8)
     (i_t^2,i_{2t}^2)(i_0^2,ci_{2t}^7)(i_t^7,ci_{3t}^2)(i_0^7,i_{3t}^7)\\

   & &(i_0^3,i_{3t}^3)(i_t^3,ci_{3t}^6)(i_0^6,ci_{2t}^3)(i_t^6,i_{2t}^6)
     (i_0^4,i_{3t}^4)(i_t^4,ci_{3t}^5)(i_0^5,ci_{2t}^4)(i_t^5,i_{2t}^5),\\

(bc)^{2^{n-4}}&=&\prod_{i=0}^{t-1}(1+i,2t-i)(2t+1+i,4t-i).
\end{array}$
\end{small}
\vskip 0.2cm

Let $A=\langle a,b,c\rangle$. Clearly, $(ab)^4=1$ and $(bc)^{2^{n-3}}=1$.
Since $ci=\frac{t}{8}-i-1$, $0\leq i\leq \frac{t}{8}-2$ if and only if $1\leq ci\leq \frac{t}{8}-1$, and since $c$ interchanges $i_{jt}^8$ and $(i+1)_{jt}^1$, it also interchanges $ci_{jt}^8$ and $c(i-1)_{jt}^1$, where $c(i-1)=\frac{t}{8}-(i-1)-1$. Thus,
 \vskip 0.2cm
 \begin{small}
$\begin{array}{lcl}
c(ab)^2c&=&\prod_{i=0}^{\frac{t}{8}-1}

    (i_{0}^{1},i_{3t}^{1})(i_{0}^{8},ci_{2t}^{1})
    (i_{t}^{1},ci_{3t}^{8})(i_{t}^{8},i_{2t}^{8})

    (i_{0}^{2},i_{3t}^{2})(i_{0}^{7},ci_{2t}^{2})
    (i_{t}^{2},ci_{3t}^{7})(i_{t}^{7},i_{2t}^{7})\\

   & &(i_{t}^{3},i_{2t}^{3})(i_{0}^{3},ci_{2t}^{6})
    (i_{t}^{6},ci_{3t}^{3})(i_{0}^{6},i_{3t}^{6})

    (i_{t}^{4},i_{2t}^{4})(i_{0}^{4},ci_{2t}^{5})
    (i_{t}^{5},ci_{3t}^{4})(i_{0}^{5},i_{3t}^{5}).\\
\end{array}$
\end{small}
\vskip 0.2cm

It is clear that $(bc)^{2^{n-4}}$ interchanges $i_0^k$ and $ci_{t}^{9-k}$ as $i_0^k+ci_{t}^{9-k}=2t+1$ (note that $1 \leq i_0^k\leq t$ and $t+1 \leq ci_{t}^{9-k}\leq 2t$), and similarly $(bc)^{2^{n-4}}$ interchanges $i_{2t}^{k}$ and $ci_{3t}^{9-k}$.
Then it is easy to check  $[(ab)^2,c]=(ab)^2c(ab)^2c=(bc)^{2^{n-4}}$.
It follows that the generators $a,b,c$ of $A$ satisfy the same relations as do $\rho_0,\rho_1,\rho_2$ in $G$, and hence $A$ is isomorphic to $G$ with order $2^n$. Clearly, $A$ is transitive and $A_1$ has order at least $4$ because $a,c\in A_1$. It follows that $|A|\geq 2^n$ and hence $|G|=2^n$. On the other hand, the generators $\r_0$, $\r_1$, $\r_2$ in $L_1:=\lg \r_0, \r_1, \r_2 \ |\ \r_0^2, \r_1^2, \r_2^2, (\r_0\r_1)^{4}, (\r_1\r_2)^{2}, (\r_0\r_2)^2\rg$ (defined in Proposition~\ref{degenerate}) satisfy all relations in $G$. This implies that $o(\r_0\r_1)=4$ in $G$, and by Proposition~\ref{stringC}, $(G,\{\r_0,\r_1,\r_2\})$ is a string C-group.

Now we prove the necessity. Let $(G,\{\r_0,\r_1,\r_2\})$ be a string C-group of rank three with type $\{4, 2^{n-3}\}$ and $|G|=2^n$. Then each of $\r_0,\r_1$ and $\r_2$ has order $2$, and we further have $o(\r_0\r_1)=4$, $o(\r_0\r_2)=2$ and $o(\r_1\r_2)=2^{n-3}$. To finish the proof, we aim to show that $G\cong G_1$ or $G_2$. Since both $G_1$ and $G_2$ are $C$-groups of order $2^n$ of type $\{4, 2^{n-3}\}$, it suffices to show that, in $G$, $[(\r_0\r_1)^2, \r_2] =1$ or $[(\r_0\r_1)^2, \r_2](\r_1\r_2)^{2^{n-4}} =1$, which will be done by induction on $n$. This can easily be checked to be true for $n = 10$ by using the computational algebra package {\sc Magma}~\cite{BCP97}.

Assume $n\geq 11$. Take $N = \lg (\r_1\r_2)^{2^{n-4}} \rg$. By Lemma~\ref{quotient}, we have $N \unlhd G$ and
($\olg=G/N, \{\overline{\r_0}, \overline{\r_{1}}, \overline{\r_{2}}\})$ (with $\overline{\r_i} = N\r_i$)
is a string C-group of rank three of type $\{4, 2^{n-4}\}$. Since $|\olg| = 2^{n-1}$, by induction hypothesis we may assume $\olg=\olg_1$ or $\olg_2$, where
\begin{small}
\begin{itemize}
  \item [$\olg_1$]$=\lg \overline{\r_0}, \overline{\r_1}, \overline{\r_2} \ |\ \overline{\r_0}^2, \overline{\r_1}^2, \overline{\r_2}^2, (\overline{\r_0}\overline{\r_1})^{2^{2}}, (\overline{\r_1}\overline{\r_2})^{2^{n-4}}, (\overline{\r_0}\overline{\r_2})^2, [(\overline{\r_0}\overline{\r_1})^2, \overline{\r_2}]\rg$,
  \item [$\olg_2$]$=\lg \overline{\r_0}, \overline{\r_1}, \overline{\r_2} \ |\ \overline{\r_0}^2, \overline{\r_1}^2, \overline{\r_2}^2, (\overline{\r_0}\overline{\r_1})^{2^{2}}, (\overline{\r_1}\overline{\r_2})^{2^{n-4}}, (\overline{\r_0}\overline{\r_2})^2, [(\overline{\r_0}\overline{\r_1})^2, \overline{\r_2}](\overline{\r_1}\overline{\r_2})^{2^{n-5}}\rg$.
\end{itemize}
\end{small}
Suppose $\olg=\olg_2$. Since $N=\langle (\r_1\r_2)^{2^{n-4}}\rangle\cong\mz_2$, we have $[(\r_0\r_1)^2, \r_2](\r_1\r_2)^{2^{n-5}}=1$ or $(\r_1\r_2)^{2^{n-4}}$, implying $[(\r_0\r_1)^2, \r_2]=(\r_1\r_2)^{\d\cdot 2^{n-5}}$, where $\d=1$ or $-1$. Since  $((\r_0\r_1)^2)^{\r_0}=(\r_0\r_1)^{-2}=(\r_0\r_1)^2$ and $\r_0\r_2=\r_2\r_0$, we have $[(\r_0\r_1)^2, \r_2]^{\r_0}=[(\r_0\r_1)^2, \r_2]$, and hence $[\r_0, (\r_1\r_2)^{\d\cdot 2^{n-5}}]=1$. By  Proposition~\ref{commutator},   $1=[(\r_0\r_1)^4,\r_2]=[(\r_0\r_1)^2,\r_2]^{(\r_0\r_1)^2}[(\r_0\r_1)^2,\r_2]=((\r_1\r_2)^{\d\cdot 2^{n-5}})^{(\r_0\r_1)^2}(\r_1\r_2)^{\d\cdot 2^{n-5}}=(\r_1\r_2)^{\d\cdot 2^{n-4}}$, which is impossible because  $o(\r_1\r_2)=2^{n-3}$.

Thus, $\olg=\olg_1$. Since $N=\langle (\r_1\r_2)^{2^{n-4}}\rangle\cong\mz_2$, we have $[(\r_0\r_1)^2, \r_2] = 1$ or $(\r_1\r_2)^{2^{n-4}}$. For the latter, $[(\r_0\r_1)^2, \r_2](\r_1\r_2)^{2^{n-4}}=(\r_1\r_2)^{2^{n-3}}=1$. It follows that $G\cong G_1$ or $G_2$.
\hfill\qed

\f{\bf Proof of Theorem~\ref{maintheorem}(2):} For the sufficiency, we need to show that both $G_3$ and $G_4$ are string $C$-group of order $2^n$, where $n \geq 10$ and
\vskip 0.1cm
\begin{small}
$\begin{array}{rl}
  G_3=& \lg \r_0, \r_1, \r_2 \ |\ \r_0^2, \r_1^2, \r_2^2, (\r_0\r_1)^{2^{2}}, (\r_1\r_2)^{2^{n-4}}, (\r_0\r_2)^2, [(\r_0\r_1)^2,(\r_1\r_2)^2] \rg,  \\
  G_4=& \lg \r_0, \r_1, \r_2 \ |\ \r_0^2, \r_1^2, \r_2^2, (\r_0\r_1)^{2^{2}}, (\r_1\r_2)^{2^{n-4}}, (\r_0\r_2)^2, [(\r_0\r_1)^2,(\r_1\r_2)^2](\r_1\r_2)^{2^{n-5}} \rg.
\end{array}$
\end{small}\\
By Corollary~\ref{cor3.1}, we only need to show that $G_4$ is a string C-group.

Let $G=G_4$. Then $[(\r_0\r_1)^2,(\r_1\r_2)^2]=(\r_1\r_2)^{2^{n-5}}$. Noting $[(\r_0\r_1)^2,(\r_1\r_2)^2(\r_1\r_2)^2]=[(\r_0\r_1)^2,(\r_1\r_2)^2][(\r_0\r_1)^2,(\r_1\r_2)^2]^{(\r_1\r_2)^2}
=(\r_1\r_2)^{2^{n-4}}=1$, we have $[(\r_0\r_1)^2,(\r_1\r_2)^{2^{n-5}}]=1$, which implies $[\r_1,((\r_1\r_2)^{2^{n-5}})^{\r_0}]=1$ as $[\r_1,(\r_1\r_2)^{2^{n-5}}]=1$. Thus, $\r_1$ fixes $K=\lg (\r_1\r_2)^{2^{n-5}}, ((\r_1\r_2)^{2^{n-5}})^{\r_0}\rg$. Clearly,  $K^{\r_0}=K^{\r_2}=K$, and so $K \trianglelefteq G$. The three generators $\r_0K$, $\r_1K$, $\r_2K$ in $G/K$ satisfy the same relations as $\r_0$, $\r_1$, $\r_2$ in $G_3$. In fact, $(\r_1K\r_2K)^{2^{n-5}}=K$, and hence $|G/K|\leq 2^{n-1}$ (here we need to check that $|G_3|=2^9$ for $n=9$ and this can be done using {\sc Magma}).
Furthermore, $[\r_1\r_2,((\r_1\r_2)^{2^{n-5}})^{\r_0}]=1$ as $[\r_2,((\r_1\r_2)^{2^{n-5}})^{\r_0}]=1$. It follows that $|K|\leq 4$ and $|G|\leq 2^{n+1}$.

Suppose $|G|=2^{n+1}$. Then $|G/K|=2^{n-1}$ and $|K|=4$.  It follows that $o(\r_1\r_2K)=2^{n-5}$ in $G/K$ and hence $o(\r_1\r_2)=2^{n-4}$ in $G$. Note that $2(n-4) \geq n$ as $n\geq 10$, by Proposition~\ref{core}, $\Core_{G}(\lg \r_1\r_2\rg)>1$, so $\lg (\r_1\r_2)^{2^{n-5}}\rg\unlhd G$. It follows that $|K|=2$, a contradiction. Thus $|G|\leq 2^n$.

\begin{figure}[H]
  \centering
  \includegraphics[width=14cm]{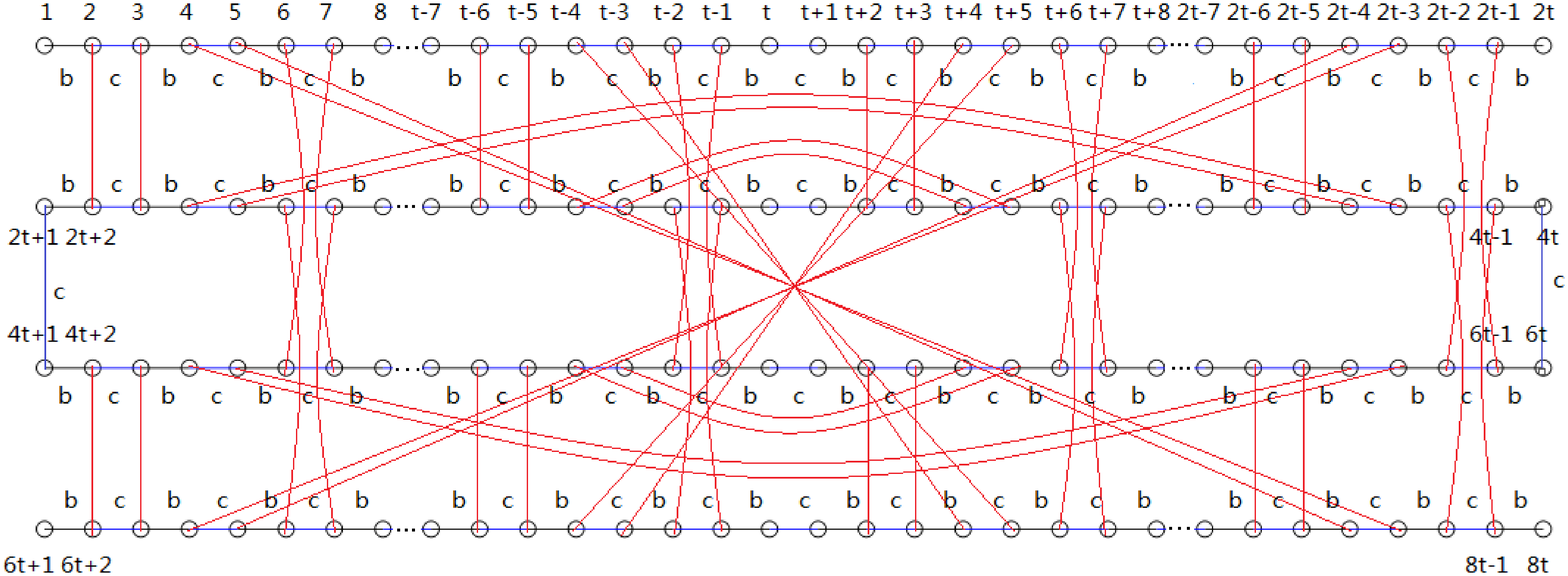}\\
  \caption{A permutation representation graph corresponding to $G_4$}\label{4}
\end{figure}

We give in Figure~\ref{4} a permutation representation graph of $G$. In this case, $bc$ consists of two paths of length $2t$ and a circle of length $4t$ with alternating labels $b$ and $c$, where $t=2^{n-5}$.

Write $i_{jt}^k=jt+8i+k$ and  ${ci}_{jt}^k=jt+8(\frac{t}{8}-i-1)+k$ for $0 \leq i \leq \frac{t}{8}-1$, $0\leq j \leq 7$ and $1\leq k \leq 8$. The permutations $a, b, c$ on the set $\{1,2,\cdots, 2^{n-2}\}$ are as follows:

\begin{small}
$\begin{array}{rl}
  a=&\prod_{i=0}^{\frac{t}{8}-1}

     (i_0^2,i_{2t}^2)(i_t^2,i_{3t}^2)(i_{4t}^2,i_{6t}^2)(i_{5t}^2,i_{7t}^2)
     (i_0^3,i_{2t}^3)(i_t^3,i_{3t}^3)(i_{4t}^3,i_{6t}^3)(i_{5t}^3,i_{7t}^3)\\

    &(i_0^4,ci_{7t}^5)(i_0^5,ci_{7t}^4)(i_{2t}^4,ci_{3t}^5)(i_{2t}^5,ci_{3t}^4)
     (i_{6t}^4,ci_{t}^{5})(i_{6t}^5,ci_{t}^{4})(i_{4t}^4,ci_{5t}^5)(i_{4t}^5,ci_{5t}^4)\\

    &(i_0^6,i_{4t}^6)(i_t^6,i_{5t}^6)(i_{2t}^6,i_{6t}^6)(i_{3t}^6,i_{7t}^6)
     (i_0^7,i_{4t}^7)(i_t^7,i_{5t}^7)(i_{2t}^7,i_{6t}^7)(i_{3t}^7,i_{7t}^7),\\
  b=&\prod_{j=0}^7\prod_{i=0}^{\frac{t}{8}-1}(i_{jt}^1,i_{jt}^2)(i_{jt}^3,i_{jt}^4)
  (i_{jt}^5,i_{jt}^6)(i_{jt}^7,i_{jt}^8),          \\
  \end{array}$
\end{small} \\

 \begin{small}
 $\begin{array}{rl}
c=&(0_{0}^{1}) (0_{6t}^{1})((\frac{t}{8}-1)_{t}^{8})((\frac{t}{8}-1)_{7t}^{8}) (0_{2t}^{1},0_{4t}^{1})((\frac{t}{8}-1)_{3t}^{8},(\frac{t}{8}-1)_{5t}^{8})\cdot\prod_{i=0}^{3}((\frac{t}{8}-1)_{2ti}^{8},0_{t+2ti}^{1})\\
   &\cdot\prod_{j=0}^7(\prod_{i=0}^{\frac{t}{8}-1}(i_{jt}^2,i_{jt}^3)(i_{jt}^4,i_{jt}^5)(i_{jt}^6,i_{jt}^7)\cdot \prod_{i=0}^{\frac{t}{8}-2}(i_{jt}^8,(i+1)_{jt}^1)).
\end{array}$
\end{small} \\
Here, $(i+1)_{jt}^k=jt+8(i+1)+k$ for $0 \leq i \leq \frac{t}{8}-2$. It is easy to see that $c$ fixes  $a$ under conjugacy, that is, $a^c=a$. It follows that $(ac)^2=1$.
Furthermore,

\begin{small}
$\begin{array}{lcl}
  ab&=&\prod_{i=0}^{\frac{t}{8}-1}

      (i_0^1,i_0^2,i_{2t}^1,i_{2t}^2)
      (i_t^1,i_t^2,i_{3t}^1,i_{3t}^2)
      (i_{4t}^1,i_{4t}^2,i_{6t}^1,i_{6t}^2)
      (i_{5t}^1,i_{5t}^2,i_{7t}^1,i_{7t}^2)\\

    & &(i_0^3,i_{2t}^4,ci_{3t}^6,ci_{7t}^5)
      (i_t^3,i_{3t}^4,ci_{2t}^6,ci_{6t}^5)
      (i_{2t}^3,i_0^4,ci_{7t}^6,ci_{3t}^5)
      (i_{3t}^3,i_t^4,ci_{6t}^6,ci_{2t}^5)

      \\

    & &(ci_{4t}^3,ci_{6t}^4,i_t^6,i_{5t}^5)
      (ci_{5t}^3,ci_{7t}^4,i_0^6,i_{4t}^5)
      (ci_{6t}^3,ci_{4t}^4,i_{5t}^6,i_t^5)
      (ci_{7t}^3,ci_{5t}^4,i_{4t}^6,i_0^5)\\

    & &(i_0^7,i_{4t}^8,i_{4t}^7,i_0^8)
      (i_t^7,i_{5t}^8,i_{5t}^7,i_{t}^8)
      (i_{2t}^7,i_{6t}^8,i_{6t}^7,i_{2t}^8)
      (i_{3t}^7,i_{7t}^8,i_{7t}^7,i_{3t}^8),\\

bc&=  &\prod_{i=0}^{1}(1+6ti,3+6ti,\cdots, 2t-1+6ti, 2t+6ti, 2t-2+6ti, \cdots, 2+6ti) \\
    & &(2t+1+2ti,2t+3+2ti,\cdots,4t-1+2ti,6t-2ti,6t-2-2ti,\cdots,4t+2-2ti), \\

(ab)^2&=&\prod_{i=0}^{\frac{t}{8}-1}
   (i_0^1,i_{2t}^1)(i_{t}^1,i_{3t}^1)
   (i_{4t}^1,i_{6t}^1)(i_{5t}^1,i_{7t}^1)
   (i_0^2,i_{2t}^2)(i_t^2,i_{3t}^2)
   (i_{4t}^2,i_{6t}^2)(i_{5t}^2,i_{7t}^2)\\

 & &(i_0^3,ci_{3t}^6)(i_t^3,ci_{2t}^6)
   (i_{2t}^3,ci_{7t}^6)(i_{3t}^3,ci_{6t}^6)
   (i_0^6,ci_{5t}^3)(i_t^6,ci_{4t}^3)
   (i_{4t}^6,ci_{7t}^3)(i_{5t}^6,ci_{6t}^3)\\

 & &(i_0^4,ci_{3t}^5)(i_t^4,ci_{2t}^5)
   (i_{2t}^4,ci_{7t}^5)(i_{3t}^4,ci_{6t}^5)
   (i_0^5,ci_{5t}^4)(i_t^5,ci_{4t}^4)
   (i_{4t}^5,ci_{7t}^4)(i_{5t}^5,ci_{6t}^4)\\

 & &(i_0^7,i_{4t}^7)(i_t^7,i_{5t}^7)
   (i_{2t}^7,i_{6t}^7)(i_{3t}^7,i_{7t}^7)
   (i_0^8,i_{4t}^8)(i_t^8,i_{5t}^8)
   (i_{2t}^8,i_{6t}^8)(i_{3t}^8,i_{7t}^8),\\

(bc)^{2^{n-5}}&=&\prod_{i=0}^{t-1}(1+i,2t-i)(6t+1+i,8t-i)\cdot\prod_{i=0}^{2t-1}(2t+1+i,6t-i),
\end{array}$
\end{small}
For $0\leq i\leq \frac{t}{8}-2$, $c$ interchanges $i_{jt}^8$ and $(i+1)_{jt}^1$, and also $ci_{jt}^8$ and $c(i-1)_{jt}^1$. Thus,

\vskip 0.2cm
\begin{small}
$\begin{array}{lcl}
c(ab)^2c&=&\prod_{i=0}^{\frac{t}{8}-1}

  (i_{0}^{1},i_{4t}^{1})(i_{t}^{1},i_{5t}^{1})
  (i_{6t}^{1},i_{2t}^{1})(i_{7t}^{1},i_{3t}^{1})

  (i_{0}^{2},ci_{3t}^{7})(i_{0}^{7},ci_{5t}^{2})
  (i_{t}^{2},ci_{2t}^{7})(i_{t}^{7},ci_{4t}^{2})\\

& &(i_{6t}^{2},ci_{5t}^{7})(i_{6t}^{7},ci_{3t}^{2})
  (i_{7t}^{2},ci_{4t}^{7})(i_{7t}^{7},ci_{2t}^{2})

  (i_{0}^{3},i_{2t}^{3})(i_{t}^{3},i_{3t}^{3})
  (i_{6t}^{3},i_{4t}^{3})(i_{7t}^{3},i_{5t}^{3})\\

& &(i_{0}^{4},ci_{5t}^{5})(i_{0}^{5},ci_{3t}^{4})
  (i_{t}^{4},ci_{4t}^{5})(i_{t}^{5},ci_{2t}^{4})

  (i_{6t}^{4},ci_{3t}^{5})(i_{6t}^{5},ci_{5t}^{4})
  (i_{7t}^{4},ci_{2t}^{5})(i_{7t}^{5},ci_{4t}^{4})\\

& &(i_{0}^{6},i_{4t}^{6})(i_{t}^{6},i_{5t}^{6})
  (i_{6t}^{6},i_{2t}^{6})(i_{7t}^{6},i_{3t}^{6})

  (i_{0}^{8},i_{2t}^{8})(i_{t}^{8},i_{3t}^{8})
  (i_{6t}^{8},i_{4t}^{8})(i_{7t}^{8},i_{5t}^{8}),\\

c(ab)^2cb&=&\prod_{i=0}^{\frac{t}{8}-1}

  (i_{0}^{1},i_{4t}^{2},ci_{t}^{8},ci_{3t}^{7})
  (ci_{t}^{1},ci_{5t}^{2},i_{0}^{8},i_{2t}^{7})
  (i_{2t}^{1},i_{6t}^{2},ci_{5t}^{8},ci_{7t}^{7})
  (ci_{3t}^{1},ci_{7t}^{2},i_{4t}^{8},i_{6t}^{7})\\

&  &(i_{4t}^{1},i_{0}^{2},ci_{3t}^{8},ci_{t}^{7})
  (ci_{5t}^{1},ci_{t}^{2},i_{2t}^{8},i_{0}^{7})
  (i_{6t}^{1},i_{2t}^{2},ci_{7t}^{8},ci_{5t}^{7})
  (ci_{7t}^{1},ci_{3t}^{2},i_{6t}^{8},i_{4t}^{7})\\

& &(i_{0}^{3},i_{2t}^{4},ci_{t}^{6},ci_{5t}^{5})
 (ci_{t}^{3},ci_{3t}^{4},i_{0}^{6},i_{4t}^{5})
 (i_{2t}^{3},i_{0}^{4},ci_{5t}^{6},ci_{t}^{5})
 (ci_{3t}^{3},ci_{t}^{4},i_{4t}^{6},i_{0}^{5}) \\

& &(i_{4t}^{3},i_{6t}^{4},ci_{3t}^{6},ci_{7t}^{5})
  (ci_{5t}^{3},ci_{7t}^{4},i_{2t}^{6},i_{6t}^{5})
  (i_{6t}^{3},i_{4t}^{4},ci_{7t}^{6},ci_{3t}^{5})
  (ci_{7t}^{3},ci_{5t}^{4},i_{6t}^{6},i_{2t}^{5}).\\

\end{array}$
\end{small}
\vskip 0.1cm

Let $A=\langle a,b,c\rangle$. It is clear that  $(bc)^{2^{n-5}}$ interchanges $i_0^j$ and $ci_t^{9-j}$ for each $1\leq j\leq 8$ because $i_0^j+ci_t^{9-j}=2t+1$ (note that $1\leq i_0^j\leq t$ and $t+1\leq ci_t^{9-j}\leq 2t$),
and similarly $(bc)^{2^{n-5}}$ also interchanges $i_{2t}^j$ and $ci_{5t}^{9-j}$, $i_{3t}^j$ and $ci_{4t}^{9-j}$, and $i_{6t}^j$ and $ci_{7t}^{9-j}$. Thus, $[c(ab)^2c,b]=(c(ab)^2cb)^2=(bc)^{2^{n-5}}$, and hence
$[(ab)^2,cbc]=(bc)^{2^{n-5}}$ as $[(bc)^{2^{n-5}},c]=1$.
Since $[(ab)^2,b]=1$, Proposition~\ref{commutator} implies  that $[(ab)^2,(bc)^2]=[(ab)^2,cbc][(ab)^2,b]^{cbc}=(bc)^{2^{n-5}}$. It follows that the generators $a,b,c$ of $A$ satisfy the same relations as do $\rho_0,\rho_1,\rho_2$ in $G$, and hence $A$ is isomorphic to $G$ with order $2^n$.

Again let $L_1=\lg \r_0, \r_1, \r_2 \ |\ \r_0^2, \r_1^2, \r_2^2, (\r_0\r_1)^{4}, (\r_1\r_2)^{2}, (\r_0\r_2)^2\rg$. The generators $\r_0, \r_1, \r_2$ in $L_1$ satisfy all relations in $G$. This means that $o(\r_0\r_1)=4$ in $G$, and by Proposition~\ref{stringC}, $(G,\{\r_0,\r_1,\r_2\})$ is a string C-group.

To prove the necessity, let $(G,\{\r_0,\r_1,\r_2\})$ be a string C-group of rank three with type $\{4, 2^{n-4}\}$ and $|G|=2^n$. Then $o(\r_0)=o(\r_1)=o(\r_2)=o(\r_0\r_2)=2$, $o(\r_0\r_1)=4$ and $o(\r_1\r_2)=2^{n-4}$. To finish the proof, we aim to show that $G\cong G_3$ or $G_4$. Since both $G_3$ and $G_4$ are $C$-groups of order $2^n$ and of type $\{4, 2^{n-4}\}$, it suffices to show that, in $G$, $[(\r_0\r_1)^2, (\r_1\r_2)^2] =1$ or $[(\r_0\r_1)^2, (\r_1\r_2)^2](\r_1\r_2)^{2^{n-5}} =1$, which will be done by induction on $n$.
This is true for $n = 10$ by {\sc Magma}.

Assume $n\geq 11$. Take $N = \lg (\r_1\r_2)^{2^{n-5}} \rg$. By Lemma~\ref{quotient}, we have $N \unlhd G$ and ($\olg=G/N, \{\overline{\r_0}, \overline{\r_{1}}, \overline{\r_{2}}\})$ (with $\overline{\r_i} = N\r_i$) is a string C-group of rank three of type $\{4, 2^{n-5}\}$. Since $|\olg| = 2^{n-1}$, by induction hypothesis we may assume $\olg=\olg_3$ or $\olg_4$, where
\begin{small}
\begin{itemize}
  \item [$\olg_3$] $= \lg \overline{\r_0}, \overline{\r_1}, \overline{\r_2} \ |\ \overline{\r_0}^2, \overline{\r_1}^2, \overline{\r_2}^2, (\overline{\r_0}\overline{\r_1})^{2^{2}}, (\overline{\r_1}\overline{\r_2})^{2^{n-5}}, (\overline{\r_0}\overline{\r_2})^2, [(\overline{\r_0}\overline{\r_1})^2, (\overline{\r_1}\overline{\r_2})^2]\rg$,

  \item [$\olg_4$] $= \lg  \overline{\r_0}, \overline{\r_1}, \overline{\r_2} \ |\ \overline{\r_0}^2, \overline{\r_1}^2, \overline{\r_2}^2, (\overline{\r_0}\overline{\r_1})^{2^{2}}, (\overline{\r_1}\overline{\r_2})^{2^{n-5}}, (\overline{\r_0}\overline{\r_2})^2, [(\overline{\r_0}\overline{\r_1})^2, (\overline{\r_1}\overline{\r_2})^2](\overline{\r_1}\overline{\r_2})^{2^{n-6}}\rg$.
\end{itemize}
\end{small}

Suppose $\olg=\olg_4$. Since $N=\lg (\r_1\r_2)^{2^{n-5}}\rg\cong \mathbb{Z}_2$, we have $[(\r_0\r_1)^2, (\r_1\r_2)^2](\r_1\r_2)^{2^{n-6}}=1$ or $(\r_1\r_2)^{2^{n-5}}$, which implies $[(\r_0\r_1)^2, (\r_1\r_2)^2]=(\r_1\r_2)^{\d \cdot 2^{n-6}}$, where $\d=1$ or $-1$. By Proposition~\ref{commutator}, $[(\r_0\r_1)^2, (\r_1\r_2)^4]=[(\r_0\r_1)^2, (\r_1\r_2)^2][(\r_0\r_1)^2, (\r_1\r_2)^2]^{(\r_1\r_2)^2}=(\r_1\r_2)^{\d \cdot 2^{n-5}}$, and $[(\r_0\r_1)^2, (\r_1\r_2)^8]=(\r_1\r_2)^{\d \cdot 2^{n-4}}=1$, implying $[(\r_0\r_1)^2, (\r_1\r_2)^{2^{n-6}}]=1$. Thus, $1=[(\r_0\r_1)^4,(\r_1\r_2)^2]=[(\r_0\r_1)^2,(\r_1\r_2)^2]^{(\r_0\r_1)^2}[(\r_0\r_1)^2,(\r_1\r_2)^2]=((\r_1\r_2)^{\d \cdot 2^{n-6}})^{(\r_0\r_1)^2}(\r_1\r_2)^{\d \cdot 2^{n-6}}$ $=(\r_1\r_2)^{\d \cdot 2^{n-5}}$, which is impossible because $o(\r_1\r_2)=2^{n-4}$.

Thus, $\olg=\olg_3$. In this case, $[(\r_0\r_1)^2, (\r_1\r_2)^2]=1$ or $(\r_1\r_2)^{2^{n-5}}$.
For the latter, $[(\r_0\r_1)^2, (\r_1\r_2)^2](\r_1\r_2)^{2^{n-5}}=(\r_1\r_2)^{2^{n-4}}=1$. It follows that $G\cong G_3$ or $G_4$. \hfill\qed

\f{\bf Proof of Theorem~\ref{maintheorem}(3):} Let $n\geq 10$ and let $G=G_5, G_6, G_7$ or $G_8$, where\\
\begin{small}
$\begin{array}{rl}
G_5=& \lg \r_0, \r_1, \r_2 \ |\ \r_0^2, \r_1^2, \r_2^2, (\r_0\r_1)^{2^{2}}, (\r_1\r_2)^{2^{n-5}}, (\r_0\r_2)^2, [\r_0, (\r_1\r_2)^2]^2, [\r_0, (\r_1\r_2)^4] \rg,  \\
G_6=& \lg \r_0, \r_1, \r_2 \ |\ \r_0^2, \r_1^2, \r_2^2, (\r_0\r_1)^{2^{2}}, (\r_1\r_2)^{2^{n-5}}, (\r_0\r_2)^2, [\r_0, (\r_1\r_2)^2]^2 (\r_1\r_2)^{2^{n-6}}, [\r_0, (\r_1\r_2)^4] \rg, \\
G_7=& \lg \r_0, \r_1, \r_2 \ |\ \r_0^2, \r_1^2, \r_2^2, (\r_0\r_1)^{2^{2}}, (\r_1\r_2)^{2^{n-5}}, (\r_0\r_2)^2, [\r_0, (\r_1\r_2)^2]^2, [\r_0, (\r_1\r_2)^4](\r_1\r_2)^{2^{n-6}} \rg, \\
G_8=& \lg \r_0, \r_1, \r_2 \ |\ \r_0^2, \r_1^2, \r_2^2, (\r_0\r_1)^{2^{2}}, (\r_1\r_2)^{2^{n-5}}, (\r_0\r_2)^2,    \\
        &      ~~~~~~~~~~~\quad                 [\r_0, (\r_1\r_2)^2]^2(\r_1\r_2)^{2^{n-6}}, [\r_0, (\r_1\r_2)^4)](\r_1\r_2)^{2^{n-6}} \rg.
\end{array}$
\end{small}

\f By Corollary~\ref{cor3.1}, we only need to show that $G_6, G_7$ and $G_8$ are string C-groups.

In all cases, $[\r_0, (\r_1\r_2)^4]=1$ or $(\r_1\r_2)^{2^{n-6}}$. It follows from Proposition~\ref{commutator} that $[\r_0, (\r_1\r_2)^8]=[\r_0, (\r_1\r_2)^4]([\r_0, (\r_1\r_2)^4])^{(\r_1\r_2)^4}=1$. Noting that $((\r_1\r_2)^8)^{\r_1}=(\r_1\r_2)^{-8}=((\r_1\r_2)^8)^{\r_2}$, we have $K=\lg (\r_1\r_2)^8\rg \unlhd G$. Clearly, $|K|\leq 2^{n-8}$, and the three generators $\r_0K$, $\r_1K$, $\r_2K$ in $G/K$ satisfy the same relations as $\r_0$, $\r_1$, $\r_2$ in $G_5$ when $n=8$. By {\sc Magma}, $|G_5|=2^8$ when $n=8$, and hence  $|G|=|G:K|\cdot|K| \leq 2^n$.

\medskip
\f {\bf Case 1:} $G=G_6$.

We construct a permutation representation graph of $G$, and in this graph, $bc$ consists four paths of length $2t$ and two circles of length $4t$ with alternating labels of $b$ and $c$, where $t=2^{n-6}$. We omit the drawing of the graph here because it is quite big.

Write $i_{jt}^k=jt+8i+k$ and  ${ci}_{jt}^k=jt+8(\frac{t}{8}-i-1)+k$, where $0 \leq i \leq \frac{t}{8}-1$, $1\leq k \leq 8$, $0\leq j \leq 15$.
The permutations $a, b, c$ on the set $\{1,2,\cdots, 2^{n-2}\}$ are as follows:

\vskip 0.2cm
\begin{small}
$\begin{array}{rl}
  a=&\prod_{i=0}^{\frac{t}{8}-1} (i_{0}^{2},i_{2t}^{2})(i_{t}^{2},i_{3t}^{2})(i_{4t}^{2},i_{7t}^{2})(i_{8t}^{2},i_{11t}^{2})(i_{12t}^{2},i_{14t}^{2})  (i_{13t}^{2},i_{15t}^{2})(i_{5t}^{2},ci_{7t}^{7})(i_{6t}^{2},ci_{2t}^{7})\\
  & (i_{9t}^{2},ci_{13t}^{7})(i_{10t}^{2},ci_{8t}^{7}) (i_{0}^{7},i_{4t}^{7})(i_{t}^{7},i_{5t}^{7})(i_{3t}^{7},i_{6t}^{7})(i_{9t}^{7},i_{12t}^{7})(i_{10t}^{7},i_{14t}^{7})         (i_{11t}^{7},i_{15t}^{7})\\
  & (i_{0}^{3},i_{2t}^{3})(i_{t}^{3},i_{3t}^{3})(i_{4t}^{3},i_{7t}^{3})(i_{8t}^{3},i_{11t}^{3}) (i_{12t}^{3},i_{14t}^{3})(i_{13t}^{3},i_{15t}^{3})(i_{5t}^{3},ci_{7t}^{6})(i_{6t}^{3},ci_{2t}^{6})\\
  &(i_{9t}^{3},ci_{13t}^{6})(i_{10t}^{3},ci_{8t}^{6})
  (i_{0}^{6},i_{4t}^{6})(i_{t}^{6},i_{5t}^{6})(i_{3t}^{6},i_{6t}^{6})(i_{9t}^{6},i_{12t}^{6})(i_{10t}^{6},i_{14t}^{6})                           (i_{11t}^{6},i_{15t}^{6})\\
  &(i_{6t}^{4},i_{8t}^{4})(i_{7t}^{4},i_{9t}^{4})(i_{0}^{4},ci_{15t}^{5})(i_{t}^{4},ci_{14t}^{5})(i_{2t}^{4},ci_{11t}^{5})(i_{3t}^{4},ci_{10t}^{5})(i_{4t}^{4},ci_{13t}^{5})
  (i_{5t}^{4},ci_{12t}^{5})\\
   &(i_{6t}^{5},i_{8t}^{5})(i_{7t}^{5},i_{9t}^{5})(i_{0}^{5},ci_{15t}^{4})(i_{t}^{5},ci_{14t}^{4})
   (i_{2t}^{5},ci_{11t}^{4})(i_{3t}^{5},ci_{10t}^{4})(i_{4t}^{5},ci_{13t}^{4})(i_{5t}^{5},ci_{12t}^{4})\\
  &(i_{8t}^{1},ci_{9t}^{8})(i_{9t}^{1},ci_{8t}^{8})(i_{10t}^{1},ci_{13t}^{8})(i_{11t}^{1},ci_{12t}^{8})(i_{12t}^{1},ci_{11t}^{8})(i_{13t}^{1},ci_{10t}^{8})
  (i_{14t}^{1},ci_{15t}^{8})(i_{15t}^{1},ci_{14t}^{8}),\\
\end{array}$
\end{small}

\begin{small}
$\begin{array}{rl}
b=&\prod_{j=0}^{15}\prod_{i=0}^{\frac{t}{8}-1}(i_{jt}^1,i_{jt}^2)(i_{jt}^3,i_{jt}^4)
  (i_{jt}^5,i_{jt}^6)(i_{jt}^7,i_{jt}^8),          \\

c=&\prod_{i=0}^{1}
(0_{8it}^{1})(0_{6t+8it}^{1})((\frac{t}{8}-1)_{t+8it}^{8})((\frac{t}{8}-1)_{7t+8it}^{8})(0_{2t+8it}^{1},0_{4t+8it}^{1})((\frac{t}{8}-1)_{3t+8it}^{8},(\frac{t}{8}-1)_{5t+8it}^{8})\\
    &\prod_{i=0}^{7}((\frac{t}{8}-1)_{2it}^{8},0_{(2i+1)t}^{1})
\cdot\prod_{j=0}^{15}(\prod_{i=0}^{\frac{t}{8}-1}(i_{jt}^2,i_{jt}^3)(i_{jt}^4,i_{jt}^5)(i_{jt}^6,i_{jt}^7) \cdot \prod_{i=0}^{\frac{t}{8}-2}(i_{jt}^8,(i+1)_{jt}^1)).
\end{array}$
\end{small}\\
Here, $(i+1)_{jt}^k=jt+8(i+1)+k$ for $0 \leq i \leq \frac{t}{8}-2$. For $0\leq i\leq \frac{t}{8}-2$, $c$ interchanges $i_{jt}^8$ and $(i+1)_{jt}^1$, and also $ci_{jt}^8$ and $c(i-1)_{jt}^1$. It is easy to see that $a$ is fixed under conjugacy by $c$, that is, $a^c=a$. It follows that $(ac)^2=1$. Furthermore

\vskip 0.2cm
\begin{small}
$\begin{array}{lcl}
  ab&=&\prod_{i=0}^{\frac{t}{8}-1}

          (i_{0}^{1},i_{0}^{2},i_{2t}^{1},i_{2t}^{2})(i_{t}^{1},i_{t}^{2},i_{3t}^{1},i_{3t}^{2})
          (i_{4t}^{1},i_{4t}^{2},i_{7t}^{1},i_{7t}^{2})(i_{5t}^{1},i_{5t}^{2},ci_{7t}^{8},ci_{7t}^{7})\\

         & &(i_{6t}^{1},i_{6t}^{2},ci_{2t}^{8},ci_{2t}^{7})(i_{8t}^{1},ci_{9t}^{7},ci_{12t}^{8},i_{11t}^{2})
          (i_{9t}^{1},ci_{8t}^{7},i_{10t}^{1},ci_{13t}^{7})(i_{11t}^{1},ci_{12t}^{7},ci_{9t}^{8},i_{8t}^{2})\\

         & &(i_{12t}^{1},ci_{11t}^{7},ci_{15t}^{8},i_{14t}^{2})(i_{13t}^{1},ci_{10t}^{7},ci_{14t}^{8},i_{15t}^{2})
          (i_{14t}^{1},ci_{15t}^{7},ci_{11t}^{8},i_{12t}^{2})(i_{15t}^{1},ci_{14t}^{7},ci_{10t}^{8},i_{13t}^{2})\\

         & &(i_{0}^{8},i_{0}^{7},i_{4t}^{8},i_{4t}^{7})(i_{t}^{8},i_{t}^{7},i_{5t}^{8},i_{5t}^{7})
          (i_{3t}^{8},i_{3t}^{7},i_{6t}^{8},i_{6t}^{7})(i_{8t}^{8},ci_{9t}^{2},i_{13t}^{8},ci_{10t}^{2})  \\

        & &(i_{0}^{3},i_{2t}^{4},ci_{11t}^{6},ci_{15t}^{5})(i_{t}^{3},i_{3t}^{4},ci_{10t}^{6},ci_{14t}^{5})
         (i_{2t}^{3},i_{0}^{4},ci_{15t}^{6},ci_{11t}^{5})(i_{3t}^{3},i_{t}^{4},ci_{14t}^{6},ci_{10t}^{5})\\

       & &(i_{4t}^{3},i_{7t}^{4},i_{9t}^{3},ci_{13t}^{5})(i_{5t}^{3},ci_{7t}^{5},ci_{9t}^{6},ci_{12t}^{5})
        (i_{6t}^{3},ci_{2t}^{5},i_{11t}^{3},i_{8t}^{4})(i_{7t}^{3},i_{4t}^{4},ci_{13t}^{6},i_{9t}^{4})\\

       & &(i_{8t}^{3},i_{11t}^{4},ci_{2t}^{6},i_{6t}^{4})(i_{10t}^{3},ci_{8t}^{5},ci_{6t}^{6},ci_{3t}^{5})
        (i_{12t}^{3},i_{14t}^{4},ci_{t}^{6},ci_{5t}^{5})(i_{13t}^{3},i_{15t}^{4},ci_{0}^{6},ci_{4t}^{5})\\

       & &(i_{14t}^{3},i_{12t}^{4},ci_{5t}^{6},ci_{t}^{5})(i_{15t}^{3},i_{13t}^{4},ci_{4t}^{6},ci_{0}^{5})
         (i_{5t}^{4},ci_{12t}^{6},ci_{9t}^{5},ci_{7t}^{6})(i_{10t}^{4},ci_{3t}^{6},ci_{6t}^{5},ci_{8t}^{6}),\\
 \end{array}$
\end{small}
\vskip 0.2cm
\begin{small}
$\begin{array}{lcl}

bc&=&\prod_{i=0}^{1}(1 +8ti,3 +8ti, \cdots, 2t-1 +8ti,2t +8ti,2t-2 +8ti, \cdots, 2 +8ti)\\
    & &(2t+1 +8ti,2t+3 +8ti, \cdots, 4t-1 +8ti, 6t +8ti,6t-2 +8ti, \cdots, 4t+2 +8ti)\\
    & &(4t+1 +8ti,4t+3 +8ti,\cdots, 6t-1 +8ti,4t +8ti,4t-2 +8ti, \cdots, 2t+2 +8ti)\\
    & &(6t+1 +8ti,6t+3 +8ti, \cdots, 8t-1 +8ti, 8t +8ti,8t-2 +8ti, \cdots, 6t+2 +8ti),
 \end{array}$
\end{small}
\vskip 0.2cm
\begin{small}
$\begin{array}{rl}

(bc)^{2^{n-6}}=&\prod_{i=0}^{t-1}(1+i,2t-i)(6t+1+i,8t-i)(8t+1+i,10t-i)(14t+1+i,16t-i) \\
     &\cdot\prod_{i=0}^{2t-1}(2t+1+i,6t-i)(10t+1+i,14t-i). \\
 \end{array}$
 \end{small}

The above computations imply $(ab)^4=1$ and $(bc)^{2^{n-5}}=1$. Furthermore,
\vskip 0.2cm
\begin{small}
$\begin{array}{lcl}
 (ab)^2&=&\prod_{i=0}^{\frac{t}{8}-1}
      (i_{0}^{1},i_{2t}^{1})(i_{t}^{1},i_{3t}^{1})(i_{4t}^{1},i_{7t}^{1})(i_{9t}^{1},i_{10t}^{1})
      (i_{5t}^{1},ci_{7t}^{8})(i_{6t}^{1},ci_{2t}^{8})(i_{8t}^{1},ci_{12t}^{8})(i_{11t}^{1},ci_{9t}^{8})\\

     &&(i_{12t}^{1},ci_{15t}^{8})(i_{13t}^{1},ci_{14t}^{8})(i_{14t}^{1},ci_{11t}^{8})(i_{15t}^{1},ci_{10t}^{8})
      (i_{0}^{8},i_{4t}^{8})(i_{t}^{8},i_{5t}^{8})(i_{3t}^{8},i_{6t}^{8})(i_{8t}^{8},i_{13t}^{8})\\

     &&(i_{0}^{2},i_{2t}^{2})(i_{t}^{2},i_{3t}^{2})(i_{4t}^{2},i_{7t}^{2})(i_{9t}^{2},i_{10t}^{2})
      (i_{5t}^{2},ci_{7t}^{7})(i_{6t}^{2},ci_{2t}^{7})(i_{8t}^{2},ci_{12t}^{7})(i_{11t}^{2},ci_{9t}^{7})\\

     &&(i_{12t}^{2},ci_{15t}^{7})(i_{13t}^{2},ci_{14t}^{7})(i_{14t}^{2},ci_{11t}^{7})(i_{15t}^{2},ci_{10t}^{7})
      (i_{0}^{7},i_{4t}^{7})(i_{t}^{7},i_{5t}^{7})(i_{3t}^{7},i_{6t}^{7})(i_{8t}^{7},i_{13t}^{7})\\

     &&(i_{4t}^{3},i_{9t}^{3})(i_{6t}^{3},i_{11t}^{3})(i_{0}^{3},ci_{11t}^{6})(i_{t}^{3},ci_{10t}^{6})
     (i_{2t}^{3},ci_{15t}^{6})(i_{3t}^{3},ci_{14t}^{6})(i_{5t}^{3},ci_{9t}^{6})(i_{7t}^{3},ci_{13t}^{6})\\

     &&(i_{8t}^{3},ci_{2t}^{6})(i_{10t}^{3},ci_{6t}^{6})(i_{12t}^{3},ci_{t}^{6})(i_{13t}^{3},i_{0}^{6})
      (i_{14t}^{3},ci_{5t}^{6})(i_{15t}^{3},ci_{4t}^{6})(i_{3t}^{6},i_{8t}^{6})(i_{7t}^{6},i_{12t}^{6})\\

    & &(i_{4t}^{4},i_{9t}^{4})(i_{6t}^{4},i_{11t}^{4})(i_{0}^{4},ci_{11t}^{5})(i_{t}^{4},ci_{10t}^{5})
     (i_{2t}^{4},ci_{15t}^{5})(i_{3t}^{4},ci_{14t}^{5})(i_{5t}^{4},ci_{9t}^{5})(i_{7t}^{4},ci_{13t}^{5})\\

     &&(i_{8t}^{4},ci_{2t}^{5})(i_{10t}^{4},ci_{6t}^{5})(i_{12t}^{4},ci_{t}^{5})(i_{13t}^{4},i_{0}^{5})
     (i_{14t}^{4},ci_{5t}^{5})(i_{15t}^{4},ci_{4t}^{5})(i_{3t}^{5},i_{8t}^{5})(i_{7t}^{5},i_{12t}^{5}),\\

c(ab)^2c&=&\prod_{i=0}^{\frac{t}{8}-1}

      (i_{0}^{1},i_{4t}^{1})(i_{t}^{1},i_{5t}^{1})(i_{3t}^{1},i_{6t}^{1})(i_{8t}^{1},i_{13t}^{1})
      (i_{2t}^{1},ci_{6t}^{8})(i_{7t}^{1},ci_{5t}^{8})(i_{9t}^{1},ci_{11t}^{8})(i_{10t}^{1},ci_{15t}^{8})\\

     & &(i_{11t}^{1},ci_{14t}^{8})(i_{12t}^{1},i_{8t}^{8})(i_{14t}^{1},ci_{13t}^{8})(i_{15t}^{1},ci_{12t}^{8})
      (i_{0}^{8},i_{2t}^{8})(i_{t}^{8},i_{3t}^{8})(i_{4t}^{8},i_{7t}^{8})(i_{9t}^{8},i_{10t}^{8})\\

     &&(i_{4t}^{2},i_{9t}^{2})(i_{6t}^{2},i_{11t}^{2})(i_{0}^{2},ci_{11t}^{7})(i_{t}^{2},ci_{10t}^{7})
     (i_{2t}^{2},ci_{15t}^{7})(i_{3t}^{2},ci_{14t}^{7})(i_{5t}^{2},ci_{9t}^{7})(i_{7t}^{2},ci_{13t}^{7})\\

     &&(i_{8t}^{2},ci_{2t}^{7})(i_{10t}^{2},ci_{6t}^{7})(i_{12t}^{2},ci_{t}^{7})(i_{13t}^{2},ci_{0}^{7})
      (i_{14t}^{2},ci_{5t}^{7})(i_{15t}^{2},ci_{4t}^{7})(i_{7t}^{7},i_{12t}^{7})(i_{8t}^{7},i_{3t}^{7})\\

     &&(i_{0}^{3},i_{2t}^{3})(i_{t}^{3},i_{3t}^{3})(i_{4t}^{3},i_{7t}^{3})(i_{9t}^{3},i_{10t}^{3})
     (i_{5t}^{3},ci_{7t}^{6})(i_{6t}^{3},ci_{2t}^{6})(i_{8t}^{3},ci_{12t}^{6})(i_{11t}^{3},ci_{9t}^{6})\\

    & &(i_{12t}^{3},ci_{15t}^{6})(i_{13t}^{3},ci_{14t}^{6})(i_{14t}^{3},ci_{11t}^{6})(i_{15t}^{3},ci_{10t}^{6})
      (i_{0}^{6},i_{4t}^{6})(i_{t}^{6},i_{5t}^{6})(i_{3t}^{6},i_{6t}^{6})(i_{8t}^{6},i_{13t}^{6})\\

     &&(i_{3t}^{4},i_{8t}^{4})(i_{7t}^{4},i_{12t}^{4})(i_{0}^{4},ci_{13t}^{5})(i_{t}^{4},ci_{12t}^{5})
     (i_{2t}^{4},ci_{8t}^{5})(i_{4t}^{4},ci_{15t}^{5})(i_{5t}^{4},ci_{14t}^{5})(i_{6t}^{4},ci_{10t}^{5})\\

    & &(i_{9t}^{4},ci_{5t}^{5})(i_{10t}^{4},ci_{t}^{5})(i_{11t}^{4},i_{0}^{5})(i_{13t}^{4},ci_{7t}^{5})
      (i_{14t}^{4},ci_{3t}^{5})(i_{15t}^{4},ci_{2t}^{5})(i_{4t}^{5},i_{9t}^{5})(i_{6t}^{5},i_{11t}^{5}),\\

\end{array}$
\end{small}

\begin{small}
$\begin{array}{lcl}
[(ab)^2,c]&=&\prod_{i=0}^{\frac{t}{8}-1}
      (i_{0}^{1},ci_{6t}^{8},ci_{t}^{8},i_{7t}^{1})
      (i_{t}^{1},i_{6t}^{1},ci_{0}^{8},ci_{7t}^{8})
      (i_{2t}^{1},i_{4t}^{1},ci_{5t}^{8},ci_{3t}^{8})
      (i_{3t}^{1},i_{5t}^{1},ci_{4t}^{8},ci_{2t}^{8})\\

    & &(i_{8t}^{1},i_{15t}^{1},ci_{9t}^{8},ci_{14t}^{8})
      (i_{9t}^{1},ci_{15t}^{8},ci_{8t}^{8},i_{14t}^{1})
      (i_{10t}^{1},ci_{11t}^{8},ci_{13t}^{8},i_{12t}^{1})
      (i_{11t}^{1},ci_{10t}^{8},ci_{12t}^{8},i_{13t}^{1})\\

    & &(i_{0}^{2},ci_{15t}^{7},ci_{t}^{7},i_{14t}^{2})
      (i_{t}^{2},ci_{14t}^{7},ci_{0}^{7},i_{15t}^{2})
      (i_{2t}^{2},ci_{11t}^{7},ci_{5t}^{7},i_{12t}^{2})
      (i_{3t}^{2},ci_{10t}^{7},ci_{4t}^{7},i_{13t}^{2})\\

    & &(i_{4t}^{2},ci_{13t}^{7},ci_{3t}^{7},i_{10t}^{2})
      (i_{5t}^{2},ci_{12t}^{7},ci_{2t}^{7},i_{11t}^{2})
      (i_{6t}^{2},i_{8t}^{2},ci_{7t}^{7},ci_{9t}^{7})
      (i_{7t}^{2},i_{9t}^{2},ci_{6t}^{7},ci_{8t}^{7})\\

    & &(i_{0}^{3},i_{14t}^{3},ci_{t}^{6},ci_{15t}^{6})
      (i_{t}^{3},i_{15t}^{3},ci_{0}^{6},ci_{14t}^{6})
      (i_{2t}^{3},i_{12t}^{3},ci_{5t}^{6},ci_{11t}^{6})
      (i_{3t}^{3},i_{13t}^{3},ci_{4t}^{6},ci_{10t}^{6})

      \\

    & &(i_{4t}^{3},i_{10t}^{3},ci_{3t}^{6},ci_{13t}^{6})
      (i_{5t}^{3},i_{11t}^{3},ci_{2t}^{6},ci_{12t}^{6})
     (i_{6t}^{3},ci_{9t}^{6},ci_{7t}^{6},i_{8t}^{3})
      (i_{7t}^{3},ci_{8t}^{6},ci_{6t}^{6},i_{9t}^{3})

      \\

    & &(i_{0}^{4},ci_{6t}^{5},ci_{t}^{5},i_{7t}^{4})
      (i_{t}^{4},i_{6t}^{4},ci_{0}^{5},ci_{7t}^{5})
      (i_{2t}^{4},i_{4t}^{4},ci_{5t}^{5},ci_{3t}^{5})
      (i_{3t}^{4},i_{5t}^{4},ci_{4t}^{5},ci_{2t}^{5})\\

    & &(i_{8t}^{4},i_{15t}^{4},ci_{9t}^{5},ci_{14t}^{5})
      (i_{9t}^{4},ci_{15t}^{5},ci_{8t}^{5},i_{14t}^{4})
      (i_{10t}^{4},ci_{11t}^{5},ci_{13t}^{5},i_{12t}^{4})
      (i_{11t}^{4},ci_{10t}^{5},ci_{12t}^{5},i_{13t}^{4}),\\

[(ab)^2,c]^b&=&\prod_{i=0}^{\frac{t}{8}-1}

     (i_{0}^{1},ci_{15t}^{8},ci_{t}^{8},i_{14t}^{1})
     (i_{t}^{1},ci_{14t}^{8},ci_{0}^{8},i_{15t}^{1})
     (i_{2t}^{1},ci_{11t}^{8},ci_{5t}^{8},i_{12t}^{1})
     (i_{3t}^{1},ci_{10t}^{8},ci_{4t}^{8},i_{13t}^{1})
     \\

    & &(i_{4t}^{1},ci_{13t}^{8},ci_{3t}^{8},i_{10t}^{1})
     (i_{5t}^{1},ci_{12t}^{8},ci_{2t}^{8},i_{11t}^{1})
     (i_{6t}^{1},i_{8t}^{1},ci_{7t}^{8},ci_{9t}^{8})
     (i_{7t}^{1},i_{9t}^{1},ci_{6t}^{8},ci_{8t}^{8})
     \\

    & &(i_{0}^{2},ci_{6t}^{7},ci_{t}^{7},i_{7t}^{2})
     (i_{t}^{2},i_{6t}^{2},ci_{0}^{7},ci_{7t}^{7})
     (i_{2t}^{2},i_{4t}^{2},ci_{5t}^{7},ci_{3t}^{7})
     (i_{3t}^{2},i_{5t}^{2},ci_{4t}^{7},ci_{2t}^{7})\\

    & &(i_{8t}^{2},i_{15t}^{2},ci_{9t}^{7},ci_{14t}^{7})
     (i_{9t}^{2},ci_{15t}^{7},ci_{8t}^{7},i_{14t}^{2})
     (i_{10t}^{2},ci_{11t}^{7},ci_{13t}^{7},i_{12t}^{2})
     (i_{11t}^{2},ci_{10t}^{7},ci_{12t}^{7},i_{13t}^{2})\\

    & &(i_{0}^{3},ci_{6t}^{6},ci_{t}^{6},i_{7t}^{3})
     (i_{t}^{3},i_{6t}^{3},ci_{0}^{6},ci_{7t}^{6})
     (i_{2t}^{3},i_{4t}^{3},ci_{5t}^{6},ci_{3t}^{6})
     (i_{3t}^{3},i_{5t}^{3},ci_{4t}^{6},ci_{2t}^{6})\\

    & &(i_{8t}^{3},i_{15t}^{3},ci_{9t}^{6},ci_{14t}^{6})
     (i_{9t}^{3},ci_{15t}^{6},ci_{8t}^{6},i_{14t}^{3})
     (i_{10t}^{3},ci_{11t}^{6},ci_{13t}^{6},i_{12t}^{3})
     (i_{11t}^{3},ci_{10t}^{6},ci_{12t}^{6},i_{13t}^{3})\\

    & &(i_{0}^{4},ci_{14t}^{4},ci_{t}^{5},ci_{15t}^{5})
     (i_{t}^{4},i_{15t}^{4},ci_{0}^{5},ci_{14t}^{5})
     (i_{2t}^{4},i_{12t}^{4},ci_{5t}^{5},ci_{11t}^{5})
     (i_{3t}^{4},i_{13t}^{4},ci_{4t}^{5},ci_{10t}^{5})\\

    & &(i_{4t}^{4},i_{10t}^{4},ci_{3t}^{5},ci_{13t}^{5})
     (i_{5t}^{4},i_{11t}^{4},ci_{2t}^{5},ci_{12t}^{5})
     (i_{6t}^{4},ci_{9t}^{5},ci_{7t}^{5},i_{8t}^{4})
     (i_{7t}^{4},ci_{8t}^{5},ci_{6t}^{5},i_{9t}^{4}),\\

\end{array}$
\end{small}
\vskip 0.2cm

Let $A=\lg a, b, c\rg$. Now, one may see that $[(ab)^2,c]^{bc}=[(ab)^2,c]^b$.
By Proposition~\ref{commutator}, $[a,(bc)^2]=[(ab)^2,c]^{bc}$ and $[a,(cb)^2]=[(ab)^2,c]^{b}$. It follows that $[a,(bc)^2]= [a,(cb)^2]$ and hence $[a,(bc)^4]=1$.
For $1 \leq j \leq 8$, it is clear that $(bc)^{2^{n-6}}$ interchanges $i_{0}^{j}$ and $ci_{t}^{9-j}$ as $i_0^j+ci_{t}^{9-j}=2t+1$ (note that $1\leq i_0^j\leq t$ and $t+1\leq ci_t^{9-j}\leq 2t$), and similarly $(bc)^{2^{n-6}}$ also
interchanges $i_{2t}^{j}$ and $ci_{5t}^{9-j}$, $i_{4t}^{j}$ and $ci_{3t}^{9-j}$, $i_{6t}^{j}$ and $ci_{7t}^{9-j}$, $i_{8t}^{j}$ and $ci_{9t}^{9-j}$, $i_{10t}^{j}$ and $ci_{13t}^{9-j}$, $i_{12t}^{j}$ and $ci_{11t}^{9-j}$, and $i_{14t}^{j}$ and $ci_{15t}^{9-j}$. This implies  $(bc)^{2^{n-6}}=([(ab)^2,c]^b)^2=([(ab)^2,c]^{bc})^2=[a,(bc)^2]^2$.
It follows that the generators $a,b,c$ of $A$ satisfy the same relations as do $\rho_0,\rho_1,\rho_2$ in $G$, and hence $A$ is isomorphic to $G$ with order $2^n$.

Again let $L_1=\lg \r_0, \r_1, \r_2 \ |\ \r_0^2, \r_1^2, \r_2^2, (\r_0\r_1)^{4}, (\r_1\r_2)^{2}, (\r_0\r_2)^2\rg$. The generators $\r_0, \r_1, \r_2$ in $L_1$ satisfy all relations in $G$. This means that $o(\r_0\r_1)=4$ in $G$, and by Proposition~\ref{stringC}, $(G,\{\r_0,\r_1,\r_2\})$ is a string C-group.

\medskip
\f {\bf Case 2:} $G=G_7$.

We construct a permutation representation graph of $G$, and in this graph, $bc$ consists one path of length $32t$ with alternating labels of $c$ and $b$, where $t=2^{n-5}$. Again, the graph is too big to be drawn in this paper.

Write $i_{jt}^k=jt+16i+k$, where $0 \leq i \leq \frac{t}{16}-1$, $1\leq k \leq 16$ and $0 \leq j \leq 1$. The permutations $a, b, c$ on the set $\{1,2,\cdots, 2^{n-5}\}$, are as follows:

\vskip 0.2cm
\begin{small}
$\begin{array}{rl}
  a=&\prod_{i=0}^{\frac{t}{16}-1}(i_0^5,i_{t}^5)(i_0^6,i_{t}^6)(i_0^7,i_{t}^7)(i_0^8,i_{t}^8)
  (i_0^9,i_{t}^9)(i_0^{10},i_{t}^{10})(i_0^{11},i_{t}^{11})(i_0^{12},i_{t}^{12})\\

  b=&(0_{0}^{1})(0_{t}^{1})((\frac{t}{16}-1)_{0}^{16},(\frac{t}{16}-1)_{t}^{16})\prod_{j=0}^{1}(\prod_{i=0}^{\frac{t}{16}-1}
  (i_{jt}^{2},i_{jt}^{3})(i_{jt}^{4},i_{jt}^{5})(i_{jt}^{6},i_{jt}^{7})
  (i_{jt}^{8},i_{jt}^{9})(i_{jt}^{10},i_{jt}^{11})\\

    &(i_{jt}^{12},i_{jt}^{13})(i_{jt}^{14},i_{jt}^{15})\cdot \prod_{i=0}^{\frac{t}{16}-2}(i_{jt}^{16},(i+1)_{jt}^{1})),          \\
  c=&\prod_{j=0}^{1}\prod_{i=0}^{\frac{t}{16}-1}(i_{jt}^{1},i_{jt}^{2})(i_{jt}^{3},i_{jt}^{4})
  (i_{jt}^{5},i_{jt}^{6})(i_{jt}^{7},i_{jt}^{8})
  (i_{jt}^{9},i_{jt}^{10})(i_{jt}^{11},i_{jt}^{12})
  (i_{jt}^{13},i_{jt}^{14})(i_{jt}^{15},i_{jt}^{16}).
\end{array}$
\end{small}\\
Here, $(i+1)_{jt}^k=jt+16(i+1)+k$ for $0 \leq i \leq \frac{t}{16}-2$. It is easy to see that $a$ is fixed under conjugacy of $c$, that is, $a^c=a$.
It follows that $(ac)^2=1$. Furthermore

\vskip 0.2cm
\begin{small}
$\begin{array}{lcl}
  ab&=&((\frac{t}{16}-1)_{0}^{16},(\frac{t}{16}-1)_{t}^{16}) \cdot \prod_{i=0}^{\frac{t}{16}-1}(i_0^2,i_0^3)(i_t^2,i_t^3)(i_0^4,i_0^5,i_{t}^4,i_{t}^5)(i_0^6,i_{t}^7)(i_0^7,i_{t}^6)(i_0^8,i_{t}^9)(i_0^9,i_{t}^8)\\

       &&(i_0^{10},i_{t}^{11})(i_0^{11},i_{t}^{10})(i_0^{12},i_{t}^{13},i_{t}^{12},i_0^{13})
       (i_0^{14},i_{0}^{15})(i_t^{15},i_{t}^{14})\cdot\prod_{j=0}^{1}\prod_{i=0}^{\frac{t}{16}-2}(i_{jt}^{16},(i+1)_{jt}^{1}),\\
  cb&=&(1,3,5, \cdots,t-1,2t,2t-2,\cdots,t+2,t+1,t+3,\cdots,2t-1,t,t-2,\cdots,4,2) \\
\end{array}$
\end{small}
\vskip 0.2cm

 The above computations imply $(ab)^4=1$ and $(bc)^{2^{n-5}}=1$. Moreover, we have
\vskip 0.2cm
\begin{small}
$\begin{array}{lcl}
(ab)^2&=&\prod_{i=0}^{\frac{t}{16}-1}(i_{0}^{4},i_{t}^{4})(i_{0}^{5},i_{t}^{5})
                         (i_{0}^{12},i_{t}^{12})(i_{0}^{13},i_{t}^{13}),\\

c(ab)^2c&=&\prod_{i=0}^{\frac{t}{16}-1}(i_{0}^{3},i_{t}^{3})(i_{0}^{6},i_{t}^{6})
                         (i_{0}^{11},i_{t}^{11})(i_{0}^{14},i_{t}^{14}),\\

[(ab)^2,c]&=&\prod_{i=0}^{\frac{t}{16}-1}(i_{0}^{3},i_{t}^{3})(i_{0}^{4},i_{t}^{4})
                         (i_{0}^{5},i_{t}^{5})(i_{0}^{6},i_{t}^{6})
                         (i_{0}^{11},i_{t}^{11})(i_{0}^{12},i_{t}^{12})
                         (i_{0}^{13},i_{t}^{13})(i_{0}^{14},i_{t}^{14}),\\

[(ab)^2,c]^b&=&\prod_{i=0}^{\frac{t}{16}-1}(i_0^2,i_{t}^2)(i_0^4,i_{t}^4)(i_0^5,i_{t}^5)(i_0^7,i_{t}^7)
(i_0^{10},i_{t}^{10})(i_0^{12},i_{t}^{12})(i_0^{13},i_{t}^{13})(i_0^{15},i_{t}^{15}),\\

[(ab)^2,c]^{bc}&=&\prod_{i=0}^{\frac{t}{16}-1}(i_0^1,i_{t}^1)(i_0^3,i_{t}^3)(i_0^6,i_{t}^6)(i_0^8,i_{t}^8)
(i_0^{9},i_{t}^{9})(i_0^{11},i_{t}^{11})(i_0^{14},i_{t}^{14})(i_0^{16},i_{t}^{16}).\\
\end{array}$
\end{small}

\vskip 0.2cm

Let $A= \lg a, ,b, c\rg$. Since $[a,c]=1$, by Proposition~\ref{commutator} we have $[a,(bc)^2]=[(ab)^2,c]^{bc}$, and hence $[a,(bc)^2]^2=1$.
The element $(bc)^{2^{n-6}}$ interchanges $i_{0}^{k}$ and $i_{t}^{k}$ as $i_{t}^{k}-i_{0}^{k}=t$ (note that $1 \leq i_{0}^{k}\leq t$ and $t+1 \leq i_{t}^{k} \leq 2t$), which implies $[(ab)^2,c]^{b}[(ab)^2,c]^{bc}=(bc)^{2^{n-6}}$.
Clearly, $[(ab)^2,c]^c=[(ab)^2,c]$.
Again by Proposition~\ref{commutator}, $[a,(bc)^4]=[a,(bc)^2][a,(bc)^2]^{(bc)^2}=[(ab)^2,c]^{bc}[(ab)^2,c]^{(bc)^3}=
([(ab)^2,c]^{cb}[(ab)^2,c]^{bc})^{(bc)^2}=([(ab)^2,c]^b[(ab)^2,c]^{bc})^{(bc)^2}
=(bc)^{2^{n-6}}$.
It follows that the generators $a,b,c$ of $A$ satisfy the same relations as do $\rho_0,\rho_1,\rho_2$ in $G$, and hence $A$ is a quotient group of $G$.
In particular, $o(bc)=2^{n-5}$ in $A$, and hence $o(\r_1\r_2)=2^{n-5}$ in $G$. It follows that $|G|=o(\r_1\r_2)^{8} \cdot |G/\lg(\r_1\r_2)^8\rg|=2^{n-8}\cdot 256=2^n$.

Again let $L_1=\lg \r_0, \r_1, \r_2 \ |\ \r_0^2, \r_1^2, \r_2^2, (\r_0\r_1)^{4}, (\r_1\r_2)^{2}, (\r_0\r_2)^2\rg$. The generators $\r_0, \r_1, \r_2$ in $L_1$ satisfy all relations in $G$. This means that $o(\r_0\r_1)=4$ in $G$, and by Proposition~\ref{stringC}, $(G,\{\r_0,\r_1,\r_2\})$ is a string C-group.

\medskip
\f {\bf Case 3:} $G=G_8$.

We construct a permutation representation graph of $G$. In this graph, $bc$ consists two paths of length $2t$ and a circle of length $4t$ alternating labels of $c$ and $b$, where $t=2^{n-6}$. Again, the graph is too big to be drawn in this paper.

Write $i_{jt}^k=jt+16i+k$ and  ${ci}_{jt}^k=jt+16(\frac{t}{16}-i-1)+k$, where $0 \leq i \leq \frac{t}{16}-1$, $1\leq k \leq 16$ and $0\leq j \leq 7$. The permutations $a, b, c$  on the set $\{1,2,\cdots, 2^{n-3}\}$ are as follows:

\vskip 0.2cm
\begin{small}
$\begin{array}{rl}
  a=&\prod_{i=0}^{\frac{t}{16}-1}
     (i_0^1,i_{4t}^1)(i_{t}^{1},i_{5t}^{1})(i_{2t}^{1},i_{6t}^{1})(i_{3t}^{1},i_{7t}^{1})
     (i_0^2,i_{4t}^2)(i_{t}^{2},i_{5t}^{2})(i_{2t}^{2},i_{6t}^{2})(i_{3t}^{2},i_{7t}^{2})\\

    &(i_{t}^{3},i_{2t}^{3})(i_{5t}^{3},i_{6t}^{3})(i_0^3,ci_{4t}^{14})(i_{3t}^{3},ci_{t}^{14})
     (i_{4t}^{3},ci_{6t}^{14})(i_{7t}^{3},ci_{3t}^{14})(i_0^{14},i_{5t}^{14})(i_{2t}^{14},i_{7t}^{14})\\

    &(i_{t}^{4},i_{2t}^{4})(i_{5t}^{4},i_{6t}^{4})(i_0^4,ci_{4t}^{13})(i_{3t}^{4},ci_{t}^{13})
     (i_{4t}^{4},ci_{6t}^{13})(i_{7t}^{4},ci_{3t}^{13})(i_0^{13},i_{5t}^{13})(i_{2t}^{13},i_{7t}^{13})\\

    &(i_{t}^{5},i_{4t}^{5})(i_{3t}^{5},i_{6t}^{5})(i_0^5,ci_{2t}^{12})(i_{2t}^{5},ci_{6t}^{12})
     (i_{5t}^{5},ci_{t}^{12})(i_{7t}^{5},ci_{5t}^{12})(i_0^{12},i_{3t}^{12})(i_{4t}^{12},i_{7t}^{12})\\

    &(i_{t}^{6},i_{4t}^{6})(i_{3t}^{6},i_{6t}^{6})(i_0^6,ci_{2t}^{11})(i_{2t}^{6},ci_{6t}^{11})
     (i_{5t}^{6},ci_{t}^{11})(i_{7t}^{6},ci_{5t}^{11})(i_0^{11},i_{3t}^{11})(i_{4t}^{11},i_{7t}^{11})\\

    &(i_0^7,ci_{5t}^{10})(i_{t}^{7},ci_{4t}^{10})(i_{2t}^{7},ci_{t}^{10})(i_{3t}^7,ci_0^{10})
     (i_{4t}^{7},ci_{7t}^{10})(i_{5t}^{7},ci_{6t}^{10})(i_{6t}^{7},ci_{3t}^{10})(i_{7t}^{7},ci_{2t}^{10})\\

    &(i_0^8,ci_{5t}^9)(i_{t}^{8},ci_{4t}^{9})(i_{2t}^{8},ci_{t}^{9})(i_{3t}^8,ci_0^9)
     (i_{4t}^{8},ci_{7t}^{9})(i_{5t}^{8},ci_{6t}^{9})(i_{6t}^{8},ci_{3t}^{9})(i_{7t}^{8},ci_{2t}^{9})\\

    &(i_0^{15},i_{2t}^{15})(i_{t}^{15},i_{3t}^{15})(i_{4t}^{15},i_{6t}^{15})(i_{5t}^{15},i_{7t}^{15})
     (i_0^{16},i_{2t}^{16})(i_{t}^{16},i_{3t}^{16})(i_{4t}^{16},i_{6t}^{16})(i_{5t}^{16},i_{7t}^{16}),\\

b=&(0_{0}^{1})((\frac{t}{16}-1)_{t}^{16})(0_{6t}^{1})((\frac{t}{16}-1)_{7t}^{16})
(0_{2t}^1,0_{4t}^{1})((\frac{t}{16}-1)_{3t}^{16},(\frac{t}{16}-1)_{5t}^{16})

\cdot\prod_{i=0}^{3}((\frac{t}{16}-1)_{2ti}^{16},0_{(2i+1)t}^{1})\\

  &\cdot \prod_{j=0}^{7}(\prod_{i=0}^{\frac{t}{16}-1}(i_{jt}^2,i_{jt}^3)(i_{jt}^4,i_{jt}^5)(i_{jt}^6,i_{jt}^7)
  (i_{jt}^8,i_{jt}^9)(i_{jt}^{10},i_{jt}^{11})(i_{jt}^{12},i_{jt}^{13})(i_{jt}^{14},i_{jt}^{15})\cdot\prod_{i=0}^{\frac{t}{16}-2}(i_{jt}^{16},(i+1)_{jt}^1)),\\

c=&\prod_{j=0}^{7}\prod_{i=0}^{\frac{t}{16}-1}(i_{jt}^{1},i_{jt}^{2})(i_{jt}^{3},i_{jt}^{4})(i_{jt}^{5},i_{jt}^{6})(i_{jt}^{7},i_{jt}^{8})
(i_{jt}^{9},i_{jt}^{10})(i_{jt}^{11},i_{jt}^{12})(i_{jt}^{13},i_{jt}^{14})(i_{jt}^{15},i_{jt}^{16}).\\
\end{array}$
\end{small}\\
Here, $(i+1)_{jt}^k=jt+16(i+1)+k$ for $0 \leq i \leq \frac{t}{16}-2$. For $0\leq i\leq \frac{t}{16}-2$, $b$ interchanges $i_{jt}^{16}$ and $(i+1)_{jt}^1$, and also $ci_{jt}^{16}$ and $c(i-1)_{jt}^1$. It is easy to see that $a$ is fixed under conjugacy of $c$, that is, $a^c=a$. It follows that $(ac)^2=1$. Furthermore

\vskip 0.2cm
\begin{small}
$\begin{array}{lcl}
ab&=  &(0_{0}^{1},0_{2t}^{1},0_{6t}^{1},0_{4t}^{1})(0_{t}^{1},(\frac{t}{16}-1)_{4t}^{16},0_{7t}^{1},(\frac{t}{16}-1)_{2t}^{16})\\

   &  &(0_{3t}^{1},(\frac{t}{16}-1)_{6t}^{16},0_{5t}^{1},(\frac{t}{16}-1)_{0}^{16})
     ((\frac{t}{16}-1)_{t}^{16},(\frac{t}{16}-1)_{5t}^{16},(\frac{t}{16}-1)_{7t}^{16},(\frac{t}{16}-1)_{3t}^{16}) \\

    &  &\prod_{i=0}^{\frac{t}{16}-2}

      ((i+1)_{0}^{1},i_{4t}^{16},(i+1)_{6t}^{1},i_{2t}^{16})
      ((i+1)_{t}^{1},i_{5t}^{16},(i+1)_{7t}^{1},i_{3t}^{16})\\

    &  &((i+1)_{2t}^{1},i_{6t}^{16},(i+1)_{4t}^{1},i_{0}^{16})
      ((i+1)_{3t}^{1},i_{7t}^{16},(i+1)_{5t}^{1},i_{t}^{16})\cdot\\

\end{array}$
\end{small}

\begin{small}
$\begin{array}{lcl}

    &  &\prod_{i=0}^{\frac{t}{16}-1}(i_{t}^{2},i_{5t}^{3},i_{6t}^{2},i_{2t}^{3})
      (i_{2t}^{2},i_{6t}^{3},i_{5t}^{2},i_{t}^{3})
      (i_{0}^{2},i_{4t}^{3},ci_{6t}^{15},ci_{4t}^{14})
      (i_{3t}^{2},i_{7t}^{3},ci_{3t}^{15},ci_{t}^{14})
      \\

    &  &(i_{4t}^{2},i_{0}^{3},ci_{4t}^{15},ci_{6t}^{14})
      (i_{7t}^{2},i_{3t}^{3},ci_{t}^{15},ci_{3t}^{14})
      (i_{0}^{14},i_{5t}^{15},i_{7t}^{14},i_{2t}^{15})
      (i_{2t}^{14},i_{7t}^{15},i_{5t}^{14},i_{0}^{15})\\

    &  &(i_{0}^{4},ci_{4t}^{12},ci_{7t}^{13},ci_{2t}^{12})
      (i_{t}^{4},i_{2t}^{5},ci_{6t}^{13},i_{4t}^{5})
      (i_{2t}^{4},i_{t}^{5},i_{4t}^{4},ci_{6t}^{12})
      (i_{3t}^{4},ci_{t}^{12},i_{5t}^{4},i_{6t}^{5}) \\

   &   &(i_{6t}^{4},i_{5t}^{5},ci_{t}^{13},i_{3t}^{5})
       (i_{7t}^{4},ci_{3t}^{12},ci_{0}^{13},ci_{5t}^{12})
       (i_{0}^{5},ci_{2t}^{13},ci_{7t}^{12},ci_{4t}^{13})
       (i_{7t}^{5},ci_{5t}^{13},ci_{0}^{12},ci_{3t}^{13})\\

   &   &(i_{0}^{6},ci_{2t}^{10},i_{7t}^{6},ci_{5t}^{10})
      (i_{t}^{6},i_{4t}^{7},ci_{7t}^{11},ci_{4t}^{10})
      (i_{2t}^{6},ci_{6t}^{10},i_{5t}^{6},ci_{t}^{10})
      (i_{3t}^{6},i_{6t}^{7},ci_{3t}^{11},ci_{0}^{10})\\

    &  &(i_{4t}^{6},i_{t}^{7},ci_{4t}^{11},ci_{7t}^{10})
      (i_{6t}^{6},i_{3t}^{7},ci_{0}^{11},ci_{3t}^{10})
      (i_{0}^{7},ci_{5t}^{11},i_{7t}^{7},ci_{2t}^{11})
      (i_{2t}^{7},ci_{t}^{11},i_{5t}^{7},ci_{6t}^{11})\\

    &  &(i_{0}^{8},ci_{5t}^{8},i_{6t}^{8},ci_{3t}^{8})
      (i_{t}^{8},ci_{4t}^{8},i_{7t}^{8},ci_{2t}^{8})
      (i_{0}^{9},ci_{3t}^{9},i_{6t}^{9},ci_{5t}^{9})
      (i_{t}^{9},ci_{2t}^{9},i_{7t}^{9},ci_{4t}^{9}),\\

cb&=  &\prod_{i=0}^{1}(1+6ti,3+6ti,\cdots, 2t-1+6ti, 2t+6ti, 2t-2+6ti, \cdots, 2+6ti) \\
   &  &(2t+1+2ti,2t+3+2ti,\cdots,4t-1+2ti,6t-2ti,6t-2-2ti,\cdots,4t+2-2ti), \\

(cb)^{2^{n-6}}&=&\prod_{i=0}^{t-1}(1+i,2t-i)(6t+1+i,8t-i)\cdot\prod_{i=0}^{2t-1}(2t+1+i,6t-i).
\end{array}$
\end{small}

The above computations imply $(ab)^4=1$, $(bc)^{2^{n-5}}=1$, and $(bc)^{2^{n-6}}=(cb)^{2^{n-6}}$. Furthermore,
\vskip 0.2cm
\begin{small}
$\begin{array}{lcl}
(ab)^2&=&\prod_{i=0}^{\frac{t}{16}-1}

       (i_{0}^{1},i_{6t}^{1})(i_{t}^{1},i_{7t}^{1})
       (i_{2t}^{1},i_{4t}^{1})(i_{3t}^{1},i_{5t}^{1})

       (i_{0}^{16},i_{6t}^{16})(i_{t}^{16},i_{7t}^{16})
       (i_{2t}^{16},i_{4t}^{16})(i_{3t}^{16},i_{5t}^{16})\\

      & &(i_{t}^{2},i_{6t}^{2})(i_{2t}^{2},i_{5t}^{2})
       (i_{0}^{2},ci_{6t}^{15})(i_{3t}^{2},ci_{3t}^{15})

       (i_{4t}^{2},ci_{4t}^{15})(i_{7t}^{2},ci_{t}^{15})
       (i_{0}^{15},i_{7t}^{15})(i_{2t}^{15},i_{5t}^{15})\\

       &&(i_{t}^{3},i_{6t}^{3})(i_{2t}^{3},i_{5t}^{3})
       (i_{0}^{3},ci_{6t}^{14})(i_{3t}^{3},ci_{3t}^{14})

       (i_{4t}^{3},ci_{4t}^{14})(i_{7t}^{3},ci_{t}^{14})
       (i_{0}^{14},i_{7t}^{14})(i_{2t}^{14},i_{5t}^{14})\\

      & &(i_{2t}^{4},i_{4t}^{4})(i_{3t}^{4},i_{5t}^{4})
       (i_{0}^{4},ci_{7t}^{13})(i_{t}^{4},ci_{6t}^{13})

       (i_{6t}^{4},ci_{t}^{13})(i_{7t}^{4},ci_{0}^{13})
       (i_{2t}^{13},i_{4t}^{13})(i_{3t}^{13},i_{5t}^{13})\\

      & &(i_{2t}^{5},i_{4t}^{5})(i_{3t}^{5},i_{5t}^{5})
       (i_{0}^{5},ci_{7t}^{12})(i_{t}^{5},ci_{6t}^{12})

       (i_{6t}^{5},ci_{t}^{12})(i_{7t}^{5},ci_{0}^{12})
       (i_{2t}^{12},i_{4t}^{12})(i_{3t}^{12},i_{5t}^{12})\\

      & &(i_{0}^{6},i_{7t}^{6})(i_{2t}^{6},i_{5t}^{6})
       (i_{t}^{6},ci_{7t}^{11})(i_{3t}^{6},ci_{3t}^{11})

       (i_{4t}^{6},ci_{4t}^{11})(i_{6t}^{6},ci_{0}^{11})
       (i_{t}^{11},i_{6t}^{11})(i_{2t}^{11},i_{5t}^{11})\\

      & &(i_{0}^{7},i_{7t}^{7})(i_{2t}^{7},i_{5t}^{7})
       (i_{t}^{7},ci_{7t}^{10})(i_{3t}^{7},ci_{3t}^{10})

       (i_{4t}^{7},ci_{4t}^{10})(i_{6t}^{7},ci_{0}^{10})
       (i_{t}^{10},i_{6t}^{10})(i_{2t}^{10},i_{5t}^{10})\\

      & &(i_{0}^{8},i_{6t}^{8})(i_{t}^{8},i_{7t}^{8})          (i_{2t}^{8},i_{4t}^{8})(i_{3t}^{8},i_{5t}^{8})

       (i_{0}^{9},i_{6t}^{9})(i_{t}^{9},i_{7t}^{9})          (i_{2t}^{9},i_{4t}^{9})(i_{3t}^{9},i_{5t}^{9})\\

\end{array}$

$\begin{array}{lcl}

c(ab)^2c&=&\prod_{i=0}^{\frac{t}{16}-1}

       (i_{t}^{1},i_{6t}^{1})(i_{2t}^{1},i_{5t}^{1})
       (i_{0}^{1},ci_{6t}^{16})(i_{3t}^{1},ci_{3t}^{16})

       (i_{4t}^{1},ci_{4t}^{16})(i_{7t}^{1},ci_{t}^{16})
       (i_{0}^{16},i_{7t}^{16})(i_{2t}^{16},i_{5t}^{16})\\

      & &(i_{0}^{2},i_{6t}^{2})(i_{t}^{2},i_{7t}^{2})
       (i_{2t}^{2},i_{4t}^{2})(i_{3t}^{2},i_{5t}^{2})

       (i_{0}^{15},i_{6t}^{15})(i_{t}^{15},i_{7t}^{15})
       (i_{2t}^{15},i_{4t}^{15})(i_{3t}^{15},i_{5t}^{15})\\

       &&(i_{2t}^{3},i_{4t}^{3})(i_{3t}^{3},i_{5t}^{3})
       (i_{0}^{3},ci_{7t}^{14})(i_{t}^{3},ci_{6t}^{14})

       (i_{6t}^{3},ci_{t}^{14})(i_{7t}^{3},ci_{0}^{14})
       (i_{2t}^{14},i_{4t}^{14})(i_{3t}^{14},i_{5t}^{14})\\

      & &(i_{t}^{4},i_{6t}^{4})(i_{2t}^{4},i_{5t}^{4})
       (i_{0}^{4},ci_{6t}^{13})(i_{3t}^{4},ci_{3t}^{13})

       (i_{4t}^{4},ci_{4t}^{13})(i_{7t}^{4},ci_{t}^{13})
       (i_{0}^{13},i_{7t}^{13})(i_{2t}^{13},i_{5t}^{13})\\

      & &(i_{0}^{5},i_{7t}^{5})(i_{2t}^{5},i_{5t}^{5})
       (i_{t}^{5},ci_{7t}^{12})(i_{3t}^{5},ci_{3t}^{12})

       (i_{4t}^{5},ci_{4t}^{12})(i_{6t}^{5},ci_{0}^{12})
       (i_{t}^{12},i_{6t}^{12})(i_{2t}^{12},i_{5t}^{12})\\

      & &(i_{2t}^{6},i_{4t}^{6})(i_{3t}^{6},i_{5t}^{6})
       (i_{0}^{6},ci_{7t}^{11})(i_{t}^{6},ci_{6t}^{11})

       (i_{6t}^{6},ci_{t}^{11})(i_{7t}^{6},i_{0}^{11})
       (i_{2t}^{11},i_{4t}^{11})(i_{3t}^{11},i_{5t}^{11})\\

      & &(i_{0}^{7},i_{6t}^{7})(i_{t}^{7},i_{7t}^{7})
       (i_{2t}^{7},i_{4t}^{7})(i_{3t}^{7},i_{5t}^{7})

       (i_{0}^{10},i_{6t}^{10})(i_{t}^{10},i_{7t}^{10})
       (i_{2t}^{10},i_{4t}^{10})(i_{3t}^{10},i_{5t}^{10})\\

      & &(i_{0}^{8},i_{7t}^{8})(i_{2t}^{8},i_{5t}^{8})
       (i_{t}^{8},ci_{7t}^{9})(i_{3t}^{8},ci_{3t}^{9})

       (i_{4t}^{8},ci_{4t}^{9})(i_{6t}^{8},ci_{0}^{9})
       (i_{t}^{9},i_{6t}^{9})(i_{2t}^{9},i_{5t}^{9})\\

\end{array}$

$\begin{array}{lcl}

[(ab)^2,c]&=&\prod_{i=0}^{\frac{t}{16}-1}

       (i_{0}^{1},i_{t}^{1},ci_{t}^{16},ci_{0}^{16})
       (i_{2t}^{1},ci_{4t}^{16},ci_{5t}^{16},i_{3t}^{1})
       (i_{4t}^{1},i_{5t}^{1},ci_{3t}^{16},ci_{2t}^{16})
       (i_{6t}^{1},ci_{6t}^{16},ci_{7t}^{16},i_{7t}^{1})\\

      & &(i_{0}^{2},ci_{0}^{15},ci_{t}^{15},i_{t}^{2})
       (i_{2t}^{2},i_{3t}^{2},ci_{5t}^{15},ci_{4t}^{15})
       (i_{4t}^{2},ci_{2t}^{15},ci_{3t}^{15},i_{5t}^{2})
       (i_{6t}^{2},i_{7t}^{2},ci_{7t}^{15},ci_{6t}^{15})\\

      & &(i_{0}^{3},i_{t}^{3},ci_{t}^{14},ci_{0}^{14})
       (i_{2t}^{3},i_{3t}^{3},ci_{5t}^{14},ci_{4t}^{14})
       (i_{4t}^{3},ci_{2t}^{14},ci_{3t}^{14},i_{5t}^{3})
       (i_{6t}^{3},ci_{6t}^{14},ci_{7t}^{14},i_{7t}^{3})\\

      & &(i_{0}^{4},ci_{0}^{13},ci_{t}^{13},i_{t}^{4})
       (i_{2t}^{4},ci_{4t}^{13},ci_{5t}^{13},i_{3t}^{4})
       (i_{4t}^{4},i_{5t}^{4},ci_{3t}^{13},ci_{2t}^{13})
       (i_{6t}^{4},i_{7t}^{4},ci_{7t}^{13},ci_{6t}^{13})\\

     &  &(i_{0}^{5},i_{t}^{5},ci_{t}^{12},ci_{0}^{12})
       (i_{2t}^{5},ci_{4t}^{12},ci_{5t}^{12},i_{3t}^{5})
       (i_{4t}^{5},i_{5t}^{5},ci_{3t}^{12},ci_{2t}^{12})
       (i_{6t}^{5},ci_{6t}^{12},ci_{7t}^{12},i_{7t}^{5})\\

      & &(i_{0}^{6},ci_{0}^{11},ci_{t}^{11},i_{t}^{6})
       (i_{2t}^{6},i_{3t}^{6},ci_{5t}^{11},ci_{4t}^{11})
       (i_{4t}^{6},ci_{2t}^{11},ci_{3t}^{11},i_{5t}^{6})
       (i_{6t}^{6},i_{7t}^{6},ci_{7t}^{11},ci_{6t}^{11})\\

      & &(i_{0}^{7},i_{t}^{7},ci_{t}^{10},ci_{0}^{10})
       (i_{2t}^{7},i_{3t}^{7},ci_{5t}^{10},ci_{4t}^{10})
       (i_{4t}^{7},ci_{2t}^{10},ci_{3t}^{10},i_{5t}^{7})
       (i_{6t}^{7},ci_{6t}^{10},ci_{7t}^{10},i_{7t}^{7})\\

      & &(i_{0}^{8},ci_{0}^{9},ci_{t}^{9},i_{t}^{8})
       (i_{2t}^{8},ci_{4t}^{9},ci_{5t}^{9},i_{3t}^{8})
       (i_{4t}^{8},i_{5t}^{8},ci_{3t}^{9},ci_{2t}^{9})
       (i_{6t}^{8},i_{7t}^{8},ci_{7t}^{9},ci_{6t}^{9}),\\

\end{array}$

$\begin{array}{lcl}

[(ab)^2,c]^b&=&\prod_{i=0}^{\frac{t}{16}-1}

       (i_{0}^{1},ci_{0}^{16},ci_{t}^{16},i_{t}^{1})
       (i_{2t}^{1},ci_{4t}^{16},ci_{5t}^{16},i_{3t}^{1})
       (i_{4t}^{1},i_{5t}^{1},ci_{3t}^{16},ci_{2t}^{16})
       (i_{6t}^{1},i_{7t}^{1},ci_{7t}^{16},ci_{6t}^{16})\\

      & &(i_{0}^{2},i_{t}^{2},ci_{t}^{15},ci_{0}^{15})
       (i_{2t}^{2},i_{3t}^{2},ci_{5t}^{15},ci_{4t}^{15})
       (i_{4t}^{2},ci_{2t}^{15},ci_{3t}^{15},i_{5t}^{2})
       (i_{6t}^{2},ci_{6t}^{15},ci_{7t}^{15},i_{7t}^{2})\\

      & &(i_{0}^{3},ci_{0}^{14},ci_{t}^{14},i_{t}^{3})
       (i_{2t}^{3},i_{3t}^{3},ci_{5t}^{14},ci_{4t}^{14})
       (i_{4t}^{3},ci_{2t}^{14},ci_{3t}^{14},i_{5t}^{3})
       (i_{6t}^{3},i_{7t}^{3},ci_{7t}^{14},ci_{6t}^{14})\\

      & &(i_{0}^{4},i_{t}^{4},ci_{t}^{13},ci_{0}^{13})
       (i_{2t}^{4},ci_{4t}^{13},ci_{5t}^{13},i_{3t}^{4})
       (i_{4t}^{4},i_{5t}^{4},ci_{3t}^{13},ci_{2t}^{13})
       (i_{6t}^{4},ci_{6t}^{13},ci_{7t}^{13},i_{7t}^{4})\\

      & &(i_{0}^{5},ci_{0}^{12},ci_{t}^{12},i_{t}^{5})
       (i_{2t}^{5},ci_{4t}^{12},ci_{5t}^{12},i_{3t}^{5})
       (i_{4t}^{5},i_{5t}^{5},ci_{3t}^{12},ci_{2t}^{12})
       (i_{6t}^{5},i_{7t}^{5},ci_{7t}^{12},ci_{6t}^{12})\\

      & &(i_{0}^{6},i_{t}^{6},ci_{t}^{11},ci_{0}^{11})
       (i_{2t}^{6},i_{3t}^{6},ci_{5t}^{11},ci_{4t}^{11})
       (i_{4t}^{6},ci_{2t}^{11},ci_{3t}^{11},i_{5t}^{6})
       (i_{6t}^{6},ci_{6t}^{11},ci_{7t}^{11},i_{7t}^{6})\\

      & &(i_{0}^{7},ci_{0}^{10},ci_{t}^{10},i_{t}^{7})
       (i_{2t}^{7},i_{3t}^{7},ci_{5t}^{10},ci_{4t}^{10})
       (i_{4t}^{7},ci_{2t}^{10},ci_{3t}^{10},i_{5t}^{7})
       (i_{6t}^{7},i_{7t}^{7},ci_{7t}^{10},ci_{6t}^{10})\\

      & &(i_{0}^{8},i_{t}^{8},ci_{t}^{9},ci_{0}^{9})
       (i_{2t}^{8},ci_{4t}^{9},ci_{5t}^{9},i_{3t}^{8})
       (i_{4t}^{8},i_{5t}^{8},ci_{3t}^{9},ci_{2t}^{9})
       (i_{6t}^{8},ci_{6t}^{9},ci_{7t}^{9},i_{7t}^{8}),\\
\end{array}$
\end{small}

\vskip 0.2cm

Let $A= \lg a, ,b, c\rg$. Now it is easy to see that $[(ab)^2,c]^c=[(ab)^2,c]^{-1}$ and $[(ab)^2,c]^{bc}=([(ab)^2,c]^b)^{-1}$. Every $4$-cycle in the product of distinct $4$-cycles of  $[(ab)^2,c]$ is either a $4$-cycle or the inverse of a $4$-cycle in $[(ab)^2,c]^b$, and $[(ab)^2,c][(ab)^2,c]^b$ is an involution, which fixes $2^{n-4}$ points including the point $1$. Then $[(ab)^2,c][(ab)^2,c]^b=[(ab)^2,c]^b[(ab)^2,c]$. It is clear that  $(cb)^{2^{n-6}}$ interchanges $i_0^j$ and $ci_t^{16-j}$ for each $1\leq j\leq 16$ because $i_0^j+ci_t^{16-j}=2t+1$ (note that $1\leq i_0^j\leq t$ and $t+1\leq ci_t^{16-j}\leq 2t$),
and similarly $(cb)^{2^{n-6}}$ also interchanges $i_{2t}^j$ and $ci_{5t}^{16-j}$, $i_{3t}^j$ and $ci_{4t}^{16-j}$, and $i_{6t}^j$ and $ci_{7t}^{16-j}$. It follows that $(cb)^{2^{n-6}}=[(ab)^2,c]^2=([(ab)^2,c]^b)^2$.

Since $[a,c]=1$, by Proposition~\ref{commutator} we have  $[a,(bc)^2]^2=([(ab)^2,c]^{bc})^2=([(ab)^2,c]^{b})^{-2}=
(cb)^{-2^{n-6}}=(bc)^{2^{n-6}}$ and $[a,(bc)^4]=[a,(bc)^2][a,(bc)^2]^{(bc)^2}=
([(ab)^2,c][(ab)^2,c]^{bcbc})^{bc}=([(ab)^2,c]([(ab)^2,c]^{-b})^{bc})^{bc}
=([(ab)^2,c]^2)^{bc}=((cb)^{2^{n-6}})^{bc}=(bc)^{2^{n-6}}$.  This implies that the generators $a,b,c$ of $A$ satisfy the same relations as do $\rho_0,\rho_1,\rho_2$ in $G$,
and hence $A$ is a quotient group of $G$.
In particular, $o(cb)=2^{n-5}$ in $A$, and hence $o(\r_1\r_2)=2^{n-5}$ in $G$. It follows that $|G|=o(\r_1\r_2)^{8} \cdot |G/\lg(\r_1\r_2)^8\rg|=2^{n-8}\cdot 256=2^n$.

Again let $L_1=\lg \r_0, \r_1, \r_2 \ |\ \r_0^2, \r_1^2, \r_2^2, (\r_0\r_1)^{4}, (\r_1\r_2)^{2}, (\r_0\r_2)^2\rg$. The generators $\r_0, \r_1, \r_2$ in $L_1$ satisfy all relations in $G$. This means that $o(\r_0\r_1)=4$ in $G$, and by Proposition~\ref{stringC}, $(G,\{\r_0,\r_1,\r_2\})$ is a string C-group.

Now we prove the necessity. Let $(G,\{\r_0,\r_1,\r_2\})$ be a string C-group of rank three with type $\{4, 2^{n-5}\}$ and $|G|=2^n$. Then each of $\r_0,\r_1$ and $\r_2$ has order $2$, and we further have $o(\r_0\r_1)=4$, $o(\r_0\r_2)=2$ and $o(\r_1\r_2)=2^{n-5}$. To finish the proof, we only need to prove $G\cong G_5, G_6, G_7$ or $G_8$. Since $G_5, G_6, G_7$ and $G_8$ are $C$-groups of order $2^n$ of type $\{4, 2^{n-5}\}$, it suffices to show that, in $G$, $[\r_0, (\r_2\r_1)^2]^2 =1$ or $[\r_0, (\r_2\r_1)^2]^2(\r_1\r_2)^{2^{n-6}} =1$, and $[\r_0, (\r_2\r_1)^4] =1$ or $[\r_0, (\r_2\r_1)^4](\r_1\r_2)^{2^{n-6}} =1$, which will be done by induction on $n$. This is true for $n = 10$ by {\sc Magma}.

Assume $n\geq 11$. Take $N = \lg (\r_1\r_2)^{2^{n-6}} \rg$. By Lemma~\ref{quotient}, we have $N \unlhd G$ and
($\olg=G/N, \{\overline{\r_0}, \overline{\r_{1}}, \overline{\r_{2}}\})$ (with $\overline{\r_i} = N\r_i$)
is a string C-group of rank three of type $\{4, 2^{n-6}\}$. Since $|\olg| = 2^{n-1}$, by induction hypothesis we may assume $\olg=\olg_5, \olg_6, \olg_7$ or $\olg_8$, where
\begin{small}
\begin{itemize}
   \item [$\olg_{5}$] $= \lg \overline{\r_0}, \overline{\r_1}, \overline{\r_2} \ |\ \overline{\r_0}^2, \overline{\r_1}^2, \overline{\r_2}^2, (\overline{\r_0}\overline{\r_1})^{4}, (\overline{\r_1}\overline{\r_2})^{2^{n-6}}, [\overline{\r_0}, (\overline{\r_1}\overline{\r_2})^2]^2, [\overline{\r_0}, (\overline{\r_1}\overline{\r_2})^4]\rg$,

  \item [$\olg_{6}$] $= \lg \overline{\r_0}, \overline{\r_1}, \overline{\r_2} \ |\ \overline{\r_0}^2, \overline{\r_1}^2, \overline{\r_2}^2, (\overline{\r_0}\overline{\r_1})^{4},(\overline{\r_1}\overline{\r_2})^{2^{n-6}}, [\overline{\r_0}, (\overline{\r_1}\overline{\r_2})^2]^2  (\overline{\r_1}\overline{\r_2})^{2^{n-7}}, [\overline{\r_0}, (\overline{\r_1}\overline{\r_2})^4]\rg$,

  \item [$\olg_{7}$] $= \lg \overline{\r_0}, \overline{\r_1}, \overline{\r_2} \ |\ \overline{\r_0}^2, \overline{\r_1}^2, \overline{\r_2}^2, (\overline{\r_0}\overline{\r_1})^{4}, (\overline{\r_1}\overline{\r_2})^{2^{n-6}}, [\overline{\r_0}, (\overline{\r_1}\overline{\r_1})^2]^2, [\overline{\r_0}, (\overline{\r_1}\overline{\r_2})^4](\overline{\r_1}\overline{\r_2})^{2^{n-7}}\rg$,

  \item [$\olg_{8}$] $= \lg \overline{\r_0}, \overline{\r_1}, \overline{\r_2} \ |\ \overline{\r_0}^2, \overline{\r_1}^2, \overline{\r_2}^2, (\overline{\r_0}\overline{\r_1})^{4}, (\overline{\r_1}\overline{\r_2})^{2^{n-6}}, [\overline{\r_0}, (\overline{\r_1}\overline{\r_2})^2]^2(\overline{\r_1}\overline{\r_2})^{2^{n-7}}, [\overline{\r_0}, (\overline{\r_1}\overline{\r_2})^4](\overline{\r_1}\overline{\r_2})^{2^{n-7}}\rg$.
\end{itemize}
\end{small}

Then $[\overline{\r_0}, (\overline{\r_1}\overline{\r_2})^4]=1$, or $[\overline{\r_0}, (\overline{\r_1}\overline{\r_2})^4](\overline{\r_1}\overline{\r_2})^{2^{n-7}}=1$, and since $N=\langle (\r_1\r_2)^{2^{n-6}}\rangle\cong\mz_2$, we have $[\r_0, (\r_1\r_2)^4]=(\r_1\r_2)^{\d\cdot2^{n-7}}$, where $\d=0,\pm1,2$. It follows $[\r_0, (\r_1\r_2)^8]=[\r_0, (\r_1\r_2)^4][\r_0, (\r_1\r_2)^4]^{(\r_1\r_2)^4}=(\r_1\r_2)^{\d\cdot2^{n-6}}$, and similarly $[\r_0, (\r_1\r_2)^{16}]=(\r_1\r_2)^{\d\cdot2^{n-5}}=1$. Since $n\geq 11$, we have  $[\r_0, (\r_1\r_2)^{2^{n-7}}]=1$, that is, $((\r_1\r_2)^{2^{n-7}})^{\r_0}=(\r_1\r_2)^{2^{n-7}}$.

Suppose $\olg=\olg_7$ or $\olg_8$. Then $[\overline{\r_0}, (\overline{\r_1}\overline{\r_2})^4](\overline{\r_1}\overline{\r_2})^{2^{n-7}}=1$, that is, $[\r_0, (\r_1\r_2)^4]=(\r_1\r_2)^{\g\cdot2^{n-7}}$ for $\g=\pm1$. It follows that $1=[\r_0^2, (\r_1\r_2)^4]=[\r_0, (\r_1\r_2)^4]^{\r_0}[\r_0, (\r_1\r_2)^4]=((\r_1\r_2)^{\g\cdot2^{n-7}})^{\r_0}(\r_1\r_2)^{\g\cdot2^{n-7}}=
(\r_1\r_2)^{\g\cdot2^{n-6}}$, which contradicts $o(\r_1\r_2)=2^{n-5}$.

Suppose $\olg=\olg_6$. Then $[\overline{\r_0}, (\overline{\r_1}\overline{\r_2})^2]^2(\overline{\r_1}\overline{\r_2})^{2^{n-7}}=1$, that is, $[\r_0, (\r_1\r_2)^2]^2=(\r_1\r_2)^{\g 2^{n-7}}$ for $\g=\pm1$. Recall that $((\r_1\r_2)^{\g 2^{n-7}})^{\r_0}=(\r_1\r_2)^{\g 2^{n-7}}$. On the other hand, $1=[\r_0^2, (\r_1\r_2)^2]=[\r_0, (\r_1\r_2)^2][\r_0, (\r_1\r_2)^2]^{\r_0}$, and hence $([\r_0, (\r_1\r_2)^2]^2)^{\r_0}=([\r_0, (\r_1\r_2)^2]^{\r_0})^2=([\r_0, (\r_1\r_2)^2]^2)^{-1}$, that is, $((\r_1\r_2)^{\g 2^{n-7}})^{\r_0}=(\r_1\r_2)^{-\g 2^{n-7}}$. It follows $(\r_1\r_2)^{\g 2^{n-7}}=(\r_1\r_2)^{-\g 2^{n-7}}$ and $(\r_1\r_2)^{\g 2^{n-6}}=1$, a contradiction.

Thus, $\olg=\olg_5$. Since $N=\langle (\r_1\r_2)^{2^{n-6}}\rangle\cong\mz_2$, we have $[\r_0, (\r_1\r_2)^2]^2=1$ or $(\r_1\r_2)^{2^{n-6}}$ and $[\r_0, (\r_1\r_2)^4]=1$ or $(\r_1\r_2)^{2^{n-6}}$. It follows that $G\cong G_5, G_6, G_7$ or $G_8$.  \hfill\qed

\medskip
\f {\bf Acknowledgements:} This work was supported by the National Natural Science Foundation of China (11571035, 11731002) and the 111 Project of China (B16002).


\begin{thebibliography}{99}
\bibitem{GroupBookss} Y. Berkovich.
 {\em Groups of Prime Power Order,
 vol. 1} (Walter de Gruyter, 2008).

\bibitem{BCP97}
W. Bosma, J. Cannon and C. Playoust.
 The {M}agma {A}lgebra {S}ystem. {I}:  the user language.
 {\em J. Symbolic Comput}. {\bf 24} (1997), 235--265.

\bibitem{CFLM2017}
P. J. Cameron, M. E. Fernandes, D. Leemans and M. Mixer.
 Highest rank of a polytope for $A_n$.
 {\em Proc. London Math. Soc}. {\bf115} (2017), 135--176.



\bibitem{atles1}M. Conder. Regular polytopes with up to 2000 flags. Available at \underline{\url{https://www.}}
    \underline{\url{math.auckland.ac.nz/~conder/RegularPolytopesWithFewFlags-ByOrder.txt}}.


\bibitem{MC}
M. Conder and D. Oliveros.
  The intersection condition for regular polytopes.
  {\em J. Combin. Theory Ser. A}. {\bf 120} (2013), 1291--1304.

\bibitem{SmallestPolytopes}
M. Conder.
  The smallest regular polytopes of given rank.
  {\em Adv. Math}. {\bf 236} (2013), 92--110.

\bibitem{HW}
H. S. M. Coxeter and W. O. J. Moser.
\newblock {\em Generators and Realitions for Discrete Groups}
\newblock (Springer-Verlag, 1972).

\bibitem{GD2016}
G. Cunningham and D. Pellicer.
Classification of tight regular polyhedra.
{\em J. Algebraic Combin}. {\bf 43} (2016),  665--691.

\bibitem{fl}
M. E. Fernandes and D. Leemans.
Polytopes of high rank for the symmetric groups.
{\em Adv. Math}. {\bf 228} (2011), 3207--3222.

\bibitem{flm1}
M. E. Fernandes, D. Leemans and M. Mixer.
 Polytopes of high rank for the alternating groups.
 {\em J. Combin. Theory Ser. A}. {\bf 119} (2012), 42--56.

\bibitem{flm2}
M. E. Fernandes, D. Leemans and M. Mixer.
 All alternating groups {$A_n$} with {$n\geq12$} have polytopes of
  rank {$\lfloor\frac{n-1}{2}\rfloor$}.
 {\em SIAM J. Discrete Math}. {\bf 26} (2012), 482--498.

\bibitem{sympolcorr}
M. E. Fernandes, D. Leemans and M. Mixer.
 Corrigendum to "Polytopes of high rank for the symmetric groups".
  {\em Adv. Math}. {\bf 238} (2013), 506--508.

\bibitem{flm}
M. E. Fernandes, D. Leemans and M. Mixer.
 Extension of the classification of high rank regular polytopes.
{\em Trans. Amer. Math. Soc}. {\bf 370} (2018), 8833--8857.

\bibitem{sc1024}
Y. Gomi, M. L. Loyola and M. L. A. N. De Las Pe\~{n}as.
String C-groups of order 1024,
{\em Contributions to Discrete Mathematics}.
{\bf 13} (2018), 1--22.

\bibitem{atles} M. I. Hartley.
 An atlas of small regular abstract polytopes. Available at
\underline{\url{http://www.}}
\underline{\url{abstract-polytopes.com/atlas/index.html}}.

\bibitem{GroupBooks} B. Huppert,
{\em Endliche Gruppen \uppercase\expandafter{\romannumeral1}}
 (Springer, 1967).


\bibitem{Loyola}
M. L. Loyola.
 String C-groups from groups of order $2^m$ and exponent at least $2^{m-3}$.
 Preprint (2008). \underline{\url{https://arxiv.org/abs/1607.01457v1}}.


\bibitem{GroupBook} I. M. Isaacs.
{\em  Finite Group Theory}
(American Mathematical Society, 2008).



\bibitem{Classfictionpgroup}
A. M. McKelden.
  Groups of order $2^m$ that contain cyclic subgroups of order $2^{m-3}$.
  {\em Amer. Math. Monthly}. {\bf 13} (1906), 121--136.


\bibitem{ARP}
P. McMullen and E. Schulte.
  {\em Abstract regular polytopes}
 (Cambridge University Press, 2002).

\bibitem{Pel2008}
D. Pellicer.
 CPR graphs and regular polytopes.
  {\em European J. Combin}. {\bf 29} (2008), 59--71.

\bibitem{Problem} E. Schulte and  A. I. Weiss.
 Problems on polytopes, their groups, and realizations.
  {\em Periodica Math. Hungarica}. {\bf 53} (2006), 231--255.

\end{thebibliography}
\end{document}